\input ./style/arxiv-general.cfg
\documentclass[aop,MSNbibl,seceqn,dvips]{arximspdf}
\makeatletter
   \@ifpackageloaded{graphicx}{}{\usepackage{graphicx}}
\makeatother
\doi{10.1214/14-AOP917} %kopijuoti is PTS
\volume{43}
\issue{4}
\pubyear{2015}
\firstpage{1712}
\lastpage{1730}
\makeatletter
\def\widetilde{\tilde}

\newcommand{\eqref}[1]{(\ref{#1})}
\def\R{{\mathbb R}}
\def\di{ \int}
\def\E{\mathbb{ E}}
\def\B{\mathbb{ B}}
\def\S{\mathbb{ S}}
\def\cal{\mathcal}
\newtheorem{Theorem}{Theorem}[section]
\newtheorem{theoremA}{Theorem}
\newtheorem{Lemma}{Lemma}[section]
\newproclaim{Definition}{Definition}[section]
\newproclaim{Remark}{Remark}[section]
\makeatother

\begin{document}
\begin{frontmatter}

%\dochead{}
\title{Integral identity and measure estimates for stationary
Fokker--Planck equations}
\runtitle{Integral identity and measure estimates}

\begin{aug}
\author[A]{\fnms{Wen} \snm{Huang}\thanksref{T1}\ead[label=e1]{wenh@mail.ustc.edu.cn}},
\author[B]{\fnms{Min} \snm{Ji}\corref{}\thanksref{T2}\ead[label=e2]{jimin@math.ac.cn}},
\author[C]{\fnms{Zhenxin} \snm{Liu}\thanksref{T3}\ead[label=e3]{zxliu@jlu.edu.cn}}
\and
\author[D]{\fnms{Yingfei} \snm{Yi}\thanksref{T4}\ead[label=e4]{yi@math.gatech.edu}}

\affiliation{University of Science and Technology of China,
Chinese Academy of Sciences,
Jilin University,
and
Georgia Institute of Technology and Jilin~University}
\address[A]{W. Huang\\
Wu Wen-Tsun Key\\
Laboratory of Mathematics\\
University of Science and\\
\quad Technology of China\\
Hefei 230026\\
People's Republic of China\\
\printead{e1}}
\thankstext{T1}{Supported in part by NSFC Grant
11225105, Fok Ying Tung Education Foundation and the Fundamental
Research Funds for the Central Universities WK0010000014.}

\address[B]{M. Ji\\
Academy of Mathematics and\hspace*{16pt}\\
\quad System Sciences\\
and\\
Hua Loo-Keng Key\\
Laboratory of Mathematics\\
Chinese Academy of Sciences\\
Beijing 100080\\
People's Republic of China\\
\printead{e2}}
\thankstext{T2}{Supported in part by NSFC Innovation Grant 10421101.}

\address[C]{Z. Liu\\
School of Mathematics\\
Jilin University\\
Changchun 130012\\
People's Republic of China\\
\printead{e3}}
\thankstext{T3}{Supported in part by NSFC Grant
11271151, SRF for ROCS, SEM, Chunmiao fund and the 985 program from Jilin
University.}

% Address of record for the research reported here
%\affiliation{Georgia Institute of Technology and Jilin University}
\address[D]{Y. Yi\\
School of Mathematics\\
Georgia Institute of Technology \\
Atlanta, Georgia 30332\\
USA\\
and\\
School of Mathematics\\
Jilin University\\
Changchun 130012\\
People's Republic of China\\
\printead{e4}}
\thankstext{T4}{Supported in part by NSF Grants
DMS-07-08331, DMS-11-09201 and a Scholarship from Jilin University.}
\runauthor{Huang, Ji, Liu and Yi}
%\author[A]{\fnms{} \snm{}\corref{}}
%\and
%\author[B]{\fnms{} \snm{}}
%\runauthor{}
%\affiliation{}
%\dedicated{}
%\address[A]{} %adresu isvedimo komanda gale!
%\address[B]{}
\end{aug}

% HISTORY:
\received{\smonth{12} \syear{2013}}
\revised{\smonth{1} \syear{2014}}
%\accepted{\smonth{} \syear{}}

% ABSTRACT
%
\begin{abstract}
We consider a Fokker--Planck equation in a general domain in
$\R^n$ with $L^{p}_{\mathrm{loc}}$ drift term and $W^{1,p}_{\mathrm{loc}}$
diffusion
term for any $p>n$. By deriving an integral identity, we give
several measure estimates of regular stationary measures in an
exterior domain with respect to diffusion and Lyapunov-like or
anti-Lyapunov-like functions. These estimates will be useful to
problems such as the existence and nonexistence of stationary
measures in a general domain as well as the concentration and limit
behaviors of stationary measures as diffusion vanishes.
\end{abstract}

% KEYWORDS
% Pirmas kwd is didziosios raides
%
\begin{keyword}[class=AMS]
\kwd[Primary ]{35Q84}
\kwd{60J60}
\kwd[; secondary ]{37B25}
\end{keyword}

\begin{keyword}
\kwd{Fokker--Planck equation}
\kwd{stationary measures}
\kwd{measure estimates}
\kwd{integral identity}
\kwd{level set method}
\end{keyword}
%
%\begin{keyword}[class=AMS]
%\kwd[Primary ]{}
%\kwd{}
%\kwd[; secondary ]{}
%\end{keyword}
%\begin{keyword}
%\kwd{}
%\end{keyword}

\end{frontmatter}

%\begin{appendix}
%\section{}
%\end{appendix}
%s1 #&#
\section{Introduction}\label{sec1}

Consider the stationary Fokker--Planck equation
%
%
%e1.1 #&#
\begin{equation}
\label{sfp}\qquad \cases{ %
\displaystyle { L}u(x)=:
\partial^2_{ij}\bigl(a^{ij}(x)u(x)\bigr)-
\partial_i\bigl(V^i(x)u(x)\bigr)=0, &\quad $x\in{\cal U},$
\vspace*{2pt}\cr
u(x)\ge0,&\quad $\displaystyle\int_{\cal U}u(x)\,\mathrm{ d}x=1,$}
\end{equation}
where $\cal U$ is a connected open set in $\R^n$ which can be
bounded, unbounded, or the entire space $\R^n$, $ L$ is the
Fokker--Planck operator, $A=(a^{ij})$ is an everywhere positive
semidefinite, $n\times n$-matrix valued function on $\cal U$,
called the diffusion matrix, and $V=(V^i)$ is a vector field on
$\cal U$ valued in $\R^n$, called the drift field. This equation is
in fact the one satisfied by stationary solutions of the
Fokker--Planck equation
%
%
%e1.2 #&#
\begin{equation}
\label{fp} \cases{ %
\displaystyle \frac{\partial u (x,t)}{\partial t} ={
L}u(x,t), & \quad $x\in{\cal U}, t>0,$
\vspace*{2pt}\cr
 u(x,t)\ge0, &\quad $\displaystyle\int_{\cal U}u(x,t)\,\mathrm{ d}x=1.$}
\end{equation}
In the above and also throughout the rest of the paper, we use short
notation $\partial_i ={\partial}/{\partial x_i}$,
$\partial^2_{ij}={\partial^2 }/{\partial x_i\,
\partial x_j}$,\vspace*{1pt} and we also adopt the usual summation convention on
$i,j=1,2,\ldots,n$
whenever applicable.

Following \cite{BR01,BR02,BRS}, etc., we make the following standard
hypothesis:

\begin{longlist}[(A)]
\item[(A)] $a^{ij}\in W^{1,p}_{\mathrm{loc}}(\cal U)$, $V^i\in L^{p}_{\mathrm{loc}}(\cal U)$ for all
$i,j=1,\ldots,n$, where $p>n$ is fixed.
\end{longlist}

Under the regularity condition (A), in the weakest situation one
considers measure solutions of \eqref{sfp}, called \emph{stationary
measures} of the Fokker--Planck equation~\eqref{fp}, which are Borel
probability measures $\mu$ satisfying
%
%
%e1.3 #&#
%e1.4 #&#
\begin{eqnarray}
V^i&\in& L^1_{\mathrm{loc}}({\cal U},\mu),\qquad  i=1,2,
\ldots,n\quad \mbox{and}\label
{Vmu}
\\
\di_{\cal U}{\cal L}f(x) \,\mathrm{ d}\mu(x)&=&0\qquad \mbox{for all }
f\in
C_0^\infty({\cal U}),\label{j2R}
\end{eqnarray}
where
\[
{\cal L}=a^{ij}\,\partial^2_{ij} +V^i\,
\partial_i
\]
is the adjoint Fokker--Planck operator and $C_0^\infty(\cal U)$
denotes the space of $C^\infty$ functions on $\cal U$ with compact
supports. If a stationary measure $\mu$ of \eqref{fp} is \emph{regular}
with density $u$, that is, $\mathrm{ d}\mu(x)=u(x)\,\mathrm{ d}x$ for some
$u\in C(\cal U)$, then it is clear that $u$ must be a \emph{weak stationary
solution} of \eqref{fp}, that is,
%
%
%e1.5 #&#
\begin{equation}
\label{j3} \cases{ %
\displaystyle\di_{\cal U}{\cal L}f(x)
u(x) \,\mathrm{ d}x=0,&\quad $\mbox{for all } f\in C_0^\infty(\cal U),$
\vspace*{2pt}\cr
u(x)\ge0, &\quad $\displaystyle\int_{\cal U}u(x)\,\mathrm{ d}x=1.$}
\end{equation}
In fact, under condition (A) and that $A=(a^{ij})$ is everywhere
positive definite in~$\cal U$, it follows from a regularity theorem
due to Bogachev--Krylov--R\"ockner \cite{BKR} that any stationary
measure $\mu$ of \eqref{fp} must admit a positive density $u\in
W^{1,p}_{\mathrm{loc}}(\cal U)$.

The purpose of the present paper is to provide several useful
measure estimates, in an exterior domain $\cal U\setminus K$ for a
compact subset $K$ of $\cal U$, of regular stationary measures of
\eqref{fp} with densities lying in $W^{1,p}_{\mathrm{loc}}(\cal U)$. Such
exterior estimates are evidentally important especially when $\cal
U$ is unbounded (e.g., ${\cal U}=\R^n, \R^n_+$) or $(a^{ij})$ is
degenerate on the boundary of $\cal U$. The measure estimates
contained in this paper are nontrivial because they do not follow
from the existing theory of elliptic equations even if $(a^{ij})$ is
everywhere positive definite in $\cal U$. Indeed, as to be seen in
the paper, measure estimates we give in this paper crucially rely on
an integral identity (Theorem \ref{l31}) which reveals fundamental
natures of stationary Fokker--Planck equations and enables one to
estimate the measure in a subdomain by making use of information
of noise distributions on the boundary of the domain. In fact, the
integral identity plays a similar role as the Pohozaev identity
does to semilinear elliptic equations. It is because of this
identity that our essential measure estimates can be made
regardless of the positive definiteness of $(a^{ij})$ in $\cal U$.

Our measure estimates in an exterior domain will be made with
respect to diffusions and derivatives of a Lyapunov-like or an
anti-Lyapunov-like function which is primarily a compact function in
the domain.

%
%de1.1 #&#
\begin{Definition}
A nonnegative function $U\in C(\cal U)$ is said to be a \emph{compact
function in $\cal U$} if:
\begin{longlist}[(ii)]
\item[(i)] $U(x)<\rho_M$, $x\in\cal U$; and
\item[(ii)] $\lim_{x\to
\partial{\cal U}} U(x)=\rho_M$,
\end{longlist}
where
$\rho_M= \sup_{x\in\cal U} U(x)$ is called the \emph{essential upper
bound of $U$}.
\end{Definition}

When $\cal U$ is unbounded, $\partial\cal U$ and the limit $x\to
\partial{\cal U}$ in (ii) above should be understood under the topology
which is defined through
a fixed homeomorphism between the extended
Euclidean space
$\E^n=\R^n\cup\partial\R^n$ and the closed unit ball
$\bar\B^n=\B^n\cup\partial\B^n$ in $\R^n$ which identifies $\R^n$
with $B^n$ and $\partial\R^n$ with
$\S^{n-1}$, and in particular, identifies each $x_*\in\S^{n-1}$ with
the infinity element $x_*^{\infty}\in\partial\R^n$ of the ray
through $x_*$.
Consequently, if ${\cal U}=\R^n$, then $x\to\partial\R^n$ under this
topology simply means $x\to\infty$ in the usual sense, and it is
easy to see that an unbounded, nonnegative function $U\in C(\R^n)$
is a compact function in $\R^n$ iff
%
%
%e1.6 #&#
\begin{equation}
\label{infty} \lim_{x\to\infty}U(x)= +\infty.
\end{equation}

For simplicity, we will use the same symbol $\Omega_\rho$ to denote
the $\rho$-sublevel set $\{x\in\cal U\dvtx U(x)<\rho\}$ of any compact
function $U$ on $\cal U$.

%
%de1.2 #&#
\begin{Definition}\label{Lya-An}
Let $U$ be a $C^2$ compact function in $\cal
U$.
\begin{longlist}[(ii)]
\item[(i)] $U$ is called a \emph{Lyapunov function} (resp., \emph
{anti-Lyapunov function}) in $\cal
U$ with respect to ${\cal L}$, if there is a $\rho_m\in
(0,\rho_M)$, called \emph{essential lower bound of $U$}, and a
constant $\gamma>0$, called \emph{Lyapunov constant} (resp., \emph
{anti-Lyapunov constant}) of $U$, such that
%
%
%e1.7 #&#
\begin{equation}
\label{U3} {\cal L} U(x)\leq-\gamma\qquad (\mbox{resp.}, \ge\gamma ), x\in
\tilde{\cal U}={\cal U}\setminus\bar\Omega_{\rho_m},
\end{equation}
where $\tilde{\cal U}$ is called \emph{essential domain of $U$}.

\item[(ii)] If $\gamma=0$ in \eqref{U3}, then $U$ is referred to
as a \emph{weak Lyapunov function} (resp., \emph{weak anti-Lyapunov
function}) in
$\cal U$ with respect to ${\cal L}$.
\end{longlist}
\end{Definition}

Below, for any $C^1$ compact function $U$ on $\cal U$ with
essential upper bound $\rho_M$, we let $h, H$ be two nonnegative,
locally bounded functions on $[0,\rho_M)$ such that
%
%
%e1.8 #&#
\begin{equation}
\label{Hupdown}\qquad h(\rho)\le a^{ij}(x)\,\partial_i U(x)\,
\partial_j U(x)\le H(\rho),\qquad  x\in U^{-1}(\rho), \rho\in[0,
\rho_M),
\end{equation}
where $U^{-1}(\rho)$ denotes the $\rho$-level set of $U$. For
instance, $h(\rho)$, respectively, $H(\rho)$, can be taken as the
infimum, respectively, the supremum, of $a^{ij}(x)$ $\partial_i
U(x)\,\partial_j U(x)$ on $U^{-1}(\rho)$. For simplicity, the
dependency of $h,H$ on $U$ will be made implicit.

For a regular stationary measure
of \eqref{fp} with density lying in $W^{1,p}_{\mathrm{loc}}(\cal U)$, the
result below gives some upper bound estimates of the measure
in the essential domain of a Lyapunov-like function in $\cal U$ with
respect to ${\cal L}$.

\renewcommand{\thetheoremA}{\Alph{theoremA}}
\begin{theoremA}\label{thmA}
 Assume \textup{(A)} and that there is
either a Lyapunov or a weak Lyapunov function $U$ in
$\cal U$ with respect to ${\cal L}$ with
essential lower, upper bound $\rho_m,\rho_M$, respectively. Then
the following hold for any regular stationary measure $\mu$ of \eqref{fp}
with density lying in $W^{1,p}_{\mathrm{loc}}(\cal U)$:
\begin{longlist}[(c)]
\item[(a)] If $U$ is a Lyapunov function with
Lyapunov constant $\gamma$, then for any $\rho_0\in(\rho_m,\rho_M)$,
there exists a constant
$C_{\rho_m,\rho_0}>0$ depending only on $\rho_m,\rho_0$ such that
\[
\mu({\cal U}\setminus\Omega_{\rho_0})\le\gamma^{-1}
C_{\rho_m,\rho_0} \Bigl(\sup_{(\rho_m,\rho_0)} H \Bigr) \mu(
\Omega_{\rho_0}\setminus\Omega_{\rho_m}),
\]
where $H$ is as in \eqref{Hupdown}.
\item[(b)] If, in addition, the Lyapunov function $U$ in \textup{(a)}
satisfies
%
%e1.9 #&#
\begin{equation}
\label{regular-value-1} \nabla U(x)\ne0\qquad \forall x\in U^{-1}(\rho)
\mbox{ for a.e. } \rho\in[\rho_m,\rho_M),
\end{equation}
then
\[
\mu({\cal U}\setminus\Omega_{\rho})\le\mathrm{ e}^{- \gamma\int
_{\rho_m}^{\rho} {1}/{H(t)}\,\mathrm{ d}t},\qquad \rho\in[
\rho_m,\rho_M), %
\]
where $H$ is as in \eqref{Hupdown}.
\item[(c)] If $U$ is a weak Lyapunov function such that
in \eqref{Hupdown} $h$ is positive and $H$ is continuous on $[\rho
_m,\rho_M)$, then for any $\rho_0\in
(\rho_m,\rho_M)$,
\[
\mu({\cal U}\setminus\Omega_{\rho_m})\le\mu( \Omega_{\rho_0}
\setminus\Omega_{\rho_m}) \mathrm{ e}^{\int_{\rho_0}^{\rho
_M}{1}/{\tilde H(\rho)}\,\mathrm{ d}\rho},
\]
where $\tilde H(\rho)=h(\rho)\int_{\rho_m}^\rho{1}/({H(s)})
\,\mathrm{ d}
s$, $\rho\in[\rho_m,\rho_M)$.
\end{longlist}
\end{theoremA}
%}

We note by Sard's theorem that if $U\in C^{n}(\cal U)$, then the
set of regular values of~$U$ is of full Lebesgue measure in
$[\rho_m,\rho_M)$, that is, \eqref{regular-value-1} is automatically
satisfied when $U\in C^{n}(\cal U)$.

For a regular stationary measure
of \eqref{fp} with density lying in $W^{1,p}_{\mathrm{loc}}(\cal U)$,
the result below gives some
lower bound estimates of the measure
in the essential domain of an anti-Lyapunov-like function in $\cal U$
with respect to ${\cal L}$.

\begin{theoremA}\label{thmB}
Assume \textup{(A)} and that there is
either an anti-Lyapunov or a weak anti-Lyapunov function $U$ in
$\cal U$ with respect to ${\cal L}$ with
essential lower, upper bound $\rho_m,\rho_M$, respectively. Then
the following hold for any regular stationary measure $\mu$ of \eqref{fp}
with density lying in $W^{1,p}_{\mathrm{loc}}(\cal U)$ and any $\rho_0\in
(\rho_m,\rho_M)$:
\begin{longlist}[(a)]

\item[(a)] If $U$ is an anti-Lyapunov function with anti-Lyapunov
constant $\gamma$
such that~\eqref{regular-value-1} holds, then
\[
\mu\bigl(\Omega_{\rho}\setminus\Omega^*_{\rho_m}\bigr)\geq\mu
\bigl(\Omega_{\rho_0}\setminus\Omega^*_{\rho_m}\bigr)
\mathrm{ e}^{\gamma
\int_{\rho_0}^\rho{1}/{H(t)}\,\mathrm{ d}t},\qquad \rho\in(\rho
_0,\rho_M),
\]
where $\Omega^*_{\rho_m}=\Omega_{\rho_m}\cup
U^{-1}(\rho_m)$ and $H$ is as in \eqref{Hupdown}.
\item[(b)] If $U$ is a weak anti-Lyapunov function such that $h$
in \eqref{Hupdown}
is positive and continuous on $[\rho_m,\rho_M)$, then
\[
\mu(\Omega_\rho\setminus\Omega_{\rho_m})\ge\mu(
\Omega_{\rho_0}\setminus\Omega_{\rho_m}) \mathrm{ e}^{\int_{\rho
_0}^\rho{1}/{\tilde H(t)}\,\mathrm{ d}t},\qquad \rho\in
[\rho_0,\rho_M),
\]
where $ \tilde H(\rho)=H(\rho)\int_{\rho_m}^{\rho}\frac{1}
{h(s)}\,\mathrm{ d}
s$, $\rho\in[\rho_m,\rho_M)$.
\end{longlist}
\end{theoremA}

Following the pioneering work of Has'minski{\u\i}
\cite{Ha2,Ha} for locally Lipschitz coefficients, the existence and
uniqueness of regular stationary measures of \eqref{fp} in $\R^n$
have been extensively studied when $(a^{ij})$ is everywhere
positive definite (see, e.g., \cite{ABR,ABG,Ben,NEW,BDR,BKR,BKR09,BR01,BR02,BRS2,BRS3,BRS1,BRS} and \cite{Sko,Ver87,Ver97,Ver99}). In
particular, Veretennikov \cite{Ver97} showed the existence when
$(a^{ij})$ is continuous and bounded under sup-norm, and $V$ is
measurable, locally bounded in $\R^n$ and satisfies
\[
V(x)\cdot x\le-\gamma,\qquad |x|\gg1
\]
for some positive
constant $\gamma$ depending on $(a^{ij})$. Later,
Bogachev--R\"ockner \cite{BR01} showed the existence and
uniqueness under condition (A) when there exists an unbounded
Lyapunov function in $\R^n$ with respect to ${\cal L}$ such that
\[
\lim_{x\to\infty}{\cal L}U(x)=-\infty.
\]
In this work, $(a^{ij})$ is even allowed to be degenerate in $\R^n$
for the existence of a stationary measure that is not necessarily
regular. Recently, Arapostathis--Borkar--Ghosh
\cite{ABG}, Theorem 2.6.10, showed the existence when
$(a^{ij}),(V^i)$ are locally Lipschitz and do not grow faster than
linearly at $\infty$, and there exists a so-called inf-compact
function satisfying \eqref{U3} in $\R^n$ with ``$\le$'' sign for
some $\gamma>0$. Bogachev--R\"ockner--Shaposhnikov \cite{BRS2}
proved the existence under condition (A) when there exists an
unbounded Lyapunov function $U$ in $\R^n$ with respect to ${\cal
L}$.

As shown in our work \cite{JM1}, the measure estimates contained in
Theorems \ref{thmA}, \ref{thmB} above are useful in dealing with problems of the
existence and nonexistence of stationary measures of \eqref{fp} in
a general domain $\cal U$ involving Lyapunov and weak Lyapunov
functions for the existence and anti-Lyapunov and weak anti-Lyapunov
functions for the nonexistence. Also, as explored in our works
\cite{JM3,JM4}, these estimates play important roles in
characterizing the concentration of stationary measures at both
global and local levels as well as in studying limit behaviors of a
family of stationary measures as diffusion matrices vanish. In
particular, even when we consider local concentration of stationary
measures defined in the entire space $\R^n$, the stationary measures
can be restricted to a subdomain in order to apply these estimates.
This is another motivation for us to consider these estimates in a
general domain.

This paper is organized as follows. In Section \ref{sec2}, we derive two
identities---an integral identity and a derivative formula, which
are of fundamental importance to the level set method to be
adopted in this paper. We prove Theorem \ref{thmA}(a)
in Section~\ref{sec3},
Theorem \ref{thmA}(b) and Theorem \ref{thmB}(a)
in Section \ref{sec4}, and Theorem \ref{thmA}(c) and
Theorem \ref{thmB}(b) in Section \ref{sec5}.

Throughout the rest of the paper, for simplicity, we will use the same
symbol $|\cdot|$ to denote absolute value of a number,
cardinality of a set and norm of a vector or a matrix.

%s2 #&#
\section{Ingredients of level set method}\label{sec2}

Our measure estimates will be carried out using the level set
method. In this section, we will prove two fundamental identities
involved in the level set method for conducting measure estimates of
stationary measures of \eqref{fp}. One is an integral identity which
will play a crucial role in capturing information of a weak
stationary solution of \eqref{fp} in each sublevel set of a
Lyapunov-like or an anti-Lyapunov-like function from its boundary.
The other one is a derivative formula which will be particularly
useful in the measure estimates of a stationary measure of
\eqref{fp} with respect to functions $h,H$ in \eqref{Hupdown}.

We call a bounded open set $\Omega$ in $\R^n$ a \emph{generalized
Lipschitz domain} if (i) it is a disjoint union of finitely many
Lipschitz subdomains; and (ii) intersections of boundaries among
these Lipschitz subdomains only occur at finitely many points.

%We say that $u\in C(\Omega)$ is a \emph{weak solution of \eqref{sfp}
%in $\Omega$} if it satisfies \eqref{j3} with $\Omega$ in place of
%$\cal U$.

%
%th2.1 #&#
\begin{Theorem}[(Integral identity)]\label{l31}
Assume that \textup{(A)} holds in a
domain $\Omega\subset\R^n$ and let $u\in W_{\mathrm{loc}}^{1, p}(\Omega)$
be a
weak stationary solution of \eqref{fp} in $\Omega$. Then for any
generalized Lipschitz domain
$\Omega^\prime\subset\subset\Omega$ and any function $F\in
C^2(\bar\Omega^\prime)$ with $F|_{\partial\Omega^\prime}={}$constant,
%
%
%e2.1 #&#
\begin{equation}
\label{31} \di_{\Omega^\prime}({\cal L}F)u \,\mathrm{ d}x= \di
_{\partial\Omega^\prime}
\bigl(a^{ij}\, \partial_i F \nu_j\bigr) u \,\mathrm{ d}s,
\end{equation}
where for a.e. $x\in\partial\Omega^\prime$,
$(\nu_j(x))$ denotes the
unit outward normal vector of $\partial\Omega^\prime$ at~$x$.
\end{Theorem}

\begin{pf} Let $F|_{\partial
\Omega^\prime}=c$ and $\Omega_*$ be a smooth domain
such that
$\Omega'\subset\subset\Omega_*\subset\subset\Omega$. Consider the
function
\[
\tilde F(x)=\cases{ %
F(x)-c, &\quad  $x\in
\Omega',$
\vspace*{2pt}\cr
0, &\quad  $x\in\Omega\setminus\Omega'.$}
\]
Clearly, $\tilde F\in W^{1, \infty}(\Omega)$ and $\operatorname{supp}(\tilde
F)\subset\bar\Omega'$. For any $0<h<1$, we let $\tilde F_h$ be the
regularization of $\tilde F$ in $\Omega$, that is,
\[
\tilde F_h(x)=h^{-n}\int_{\Omega}\xi
\biggl(\frac{x-y}{h} \biggr)\tilde F(y) \,\mathrm{ d}y,
\]
where the function $\xi$ is a mollifier---a nonnegative $C^\infty$
function in $\R^n$ vanishing outside of the unit ball of $\R^n$
centered at the origin and satisfying $\int_{\R^n}\xi(x)\,\mathrm{d}x=1$. Then
$\tilde F_h\in C_0^\infty(\Omega)$, $\operatorname{supp}(\tilde F_h)\subset
\bar\Omega_*$ as $0<h\ll1$, and $\tilde F_h\to\tilde F$ in $ W^{1,
q}(\Omega_*)$, as $h\to0$, for any $0<q<\infty$. Since $u$ is a
weak stationary solution of \eqref{fp} in $\Omega$,
%
%
%e2.2 #&#
\begin{equation}
\label{fh} \di_{\Omega}\bigl(a^{ij}\, \partial^2_{ij}
\tilde F_h + V^i \,\partial_i \tilde
F_h\bigr)u \,\mathrm{ d}x=0 \qquad\mbox{as } 0<h\ll1.
\end{equation}
We note that $a^{ij}, u\in W^{1,
p}_{\mathrm{loc}}(\Omega)$, $i,j=1,2,\ldots,n$. We have by passing to the
limit $h\to0$ that
%
%
%e2.3 #&#
\begin{eqnarray}
\label{25q} \di_{\Omega} u a^{ij}\, \partial^2_{ij}
\tilde F_h \,\mathrm{ d}x &=&\di_{\Omega_*} u a^{ij}\,
\partial^2_{ij} \tilde F_h \,\mathrm{ d}x=-
\di_{\Omega_*}\, \partial_j \bigl(u a^{ij}\bigr) (
\partial_i \tilde F_h) \,\mathrm{ d}x
\nonumber
\\
&\to& -\di_{\Omega_*}\,\partial_j \bigl(u a^{ij}\bigr) (
\partial_i \tilde F) \,\mathrm{ d}x =-\di_{\Omega'}\,\partial_j
\bigl(u a^{ij}\bigr) (\partial_i F) \,\mathrm{ d}x
\\
&=& \di_{\Omega'} u a^{ij}\, \partial^2_{ij} F
\,\mathrm{ d}x - \di_{\partial\Omega'} u a^{ij}\, \partial_i F
\nu_j \,\mathrm{ d}s.\nonumber
\end{eqnarray}
On the other hand, we note by the Sobolev embedding theorem that
$u\in C (\bar{\Omega}_*)$, and hence $V^iu\in L^{p}(\Omega_*)$,
$i=1,2,\ldots,n$. Thus, we can also pass to the limit $h\to0$ to
obtain
%
%
%e2.4 #&#
\begin{eqnarray}\qquad
\label{26q} \di_{\Omega}u V^i\,\partial_i \tilde
F_h \,\mathrm{ d}x&=&\di_{\Omega_*}u V^i\,\partial_i
\tilde F_h \,\mathrm{ d}x
\nonumber
\\[-8pt]
\\[-8pt]
\nonumber
&\to&\di_{\Omega_*}u V^i\,\partial_i \tilde F \,\mathrm{ d}x=
\di_{\Omega'}u V^i\,\partial_i \tilde F \,\mathrm{ d}x =
\di_{\Omega'}u V^i\,\partial_i F \,\mathrm{ d}x.
\end{eqnarray}
The theorem now follows from \eqref{fh}--(\ref{26q}).
\end{pf}

%
%re2.1 #&#
\begin{Remark} 1. We note that the theorem does not require
$(a^{ij})$ to be even positive semidefinite. It also holds for less
regular $(a^{ij})$, $(V^i)$, and $u$, as long as $a^{ij}u\in
W^{1,\alpha}_{\mathrm{loc}}(\Omega)$ and $V^iu\in L^\alpha_{\mathrm{loc}}(\Omega)$,
$\forall$ $i,j,=1,2,\ldots,n$, for some $\alpha>1$.

2. In applying the integral identity (\ref{31}), one typically
chooses $\Omega^\prime$ as a sublevel set of a Lyapunov-like or an
anti-Lyapunov-like function $U$. Of course, $\Omega^\prime$, being
such a sublevel set, need not be a generalized Lipschitz domain. As
will be seen in the next section, a technique to get around that is to
use the approximation of $U$ by Morse functions.
\end{Remark}

%
%th2.2 #&#
\begin{Theorem}[(Derivative formula)] \label{deri}
Let $\mu$ be a Borel probability measure with density
$u\in C(\cal U)$. For a compact function $U\in C^1(\cal U)$,
consider the measure function
\[
y(\rho):=\mu(\Omega_\rho)=\di_{\Omega_\rho} u \,\mathrm{ d}x,\qquad
\rho\in(0,
\rho_M), %
\]
and the open set
%
%
%e2.5 #&#
\begin{equation}
\label{regular-value} {\cal I}=:\bigl\{\rho\in(0,\rho_M)\dvtx
\nabla U(x)
\ne0, x\in U^{-1}(\rho)\bigr\},
\end{equation}
where $\rho_M$ is the essential upper bound of $U$ and $\Omega_\rho$
is the $\rho$-sublevel set of $U$ for each $\rho\in(0,\rho_M)$.
Then $y$ is of the class $C^1$ on $ \cal I$ with derivatives
%
%
%e2.6 #&#
\begin{equation}
\label{y'} y^{\prime}(\rho)=\int_{\partial\Omega_{\rho}}
\frac{u}{|\nabla U|} \,\mathrm{ d}s,\qquad \rho\in\cal I.
\end{equation}
\end{Theorem}

\begin{pf} Since $U$ is a compact function on $\cal U$, it is easy to
see that
$\partial\Omega_\rho\subset U^{-1}(\rho)$ for all $\rho\in
(0,\rho_M)$.
Let $\rho\in\cal I$. Then $\nabla U(x)\neq0$,
$x\in\partial\Omega_\rho$. Hence,
$\partial\Omega_\rho$ is
a $C^1$ hypersurface which coincides with $U^{-1}(\rho)$.

Let ${\cal T}|=\{(x,e_j)\dvtx j=1,2,\ldots,n\}$ be an orientation
preserving, orthonormal, moving frame defined over $
\partial\Omega_\rho$ such that for each $x\in\partial\Omega_\rho$,
$e_j=e_j(x)$,
$j=1,2,\ldots,n-1$, are tangent vectors, and $e_n=e_n(x)$ is the
outward unit normal vector, of $\partial\Omega_\rho$ at $x$. We
denote $\{(x,\omega^j)\dvtx j=1,2,\ldots,n\}$ as the dual frame of
$\cal T$ defined over $\partial\Omega_\rho$, that is,
$\omega^i(e_j)=\delta^i_{j}$, $i,j=1,2,\ldots,n$. Since, for each
$x\in
\partial\Omega_\rho$, $e_n={\nabla U}/{|\nabla U|}$, we have
$\omega^n={\mathrm{ d}U}/{|\nabla U|}$. Therefore,
\[
\mathrm{ d}x=\mathrm{ d}x_1\wedge\cdots\wedge\,\mathrm{ d}x_n=\mathrm{ d}s\wedge
\omega^n=\frac{1}{|\nabla
U|}\,\mathrm{ d}s \,\mathrm{ d}U, %
\]
where $\mathrm{ d}
s=\omega^1\wedge\cdots\wedge\omega^{n-1}$ is a volume form defined
on $\partial\Omega_\rho$, from which \eqref{y'} easily follows.

Continuity of $y^{\prime}(\rho)$ on $\cal I$ follows from
\eqref{y'}.
\end{pf}

%
%re2.2 #&#
\begin{Remark} In fact, the derivative formula \eqref{y'} is
known when\break $|\nabla U(x)|\ge c>0$ a.e. in $\R^n$ (see
\cite{Boga}, Proposition 5.8.34), and is already used in
\cite{BKS1}, Proposition~2, for level set estimates concerning
functions that satisfy \eqref{infty}.
\end{Remark}

%s3 #&#
\section{Proof of Theorem \texorpdfstring{\protect\ref{thmA}(a)}{A(a)}}\label{sec3}

Let $U$ be a Lyapunov function in
$\cal U$ with respect to ${\cal L}$ with
Lyapunov constant $\gamma$ and essential lower bound $\rho_m$ and
upper bound $\rho_M$
and let
$\Omega_\rho$ denote the $\rho$-sublevel set of $U$ for each
$\rho\in[\rho_m,\rho_M)$.

Given $\rho\in[\rho_m,\rho_M)$, we fix a $\rho^*\in
(\rho,\rho_M)$. Since Morse functions are dense in $C^2(\cal U)$,
there is a sequence $U_k\in C^2(\cal U)$, $k=1,2,\ldots,$ of Morse
functions such that $U_k\to U$ in
$C^2(\cal U)$, in particular, $U_k\to U$ in
$C^2(\bar\Omega_{\rho^*})$, as $k\to\infty$.
For each $k$, denote
\begin{eqnarray*}
&&\Omega^k_{\rho}=\bigl\{x\in\Omega_{\rho^*}\dvtx\mbox{
either }U_k(x)<\rho\mbox{ or } x \mbox{ is a local maximal point of }
\\
&& \hspace*{205pt}U_k \mbox{ lying in } U_k^{-1}(\rho)\bigr\}.
\end{eqnarray*}
It is obvious that $\Omega^k_{\rho}$'s are nonempty open sets
for all $k\ge1$.

%
%le3.1 #&#
\begin{Lemma}\label{Claim-1} There is a positive integer $k(\rho)$
such that
$\Omega^k_\rho\subset\subset\Omega_{\rho^\ast}$ for all $k\ge
k(\rho)$.
\end{Lemma}

\begin{pf} If this is not true, then there are sequences $k_i\to\infty$,
$x_{i}\in\partial\Omega^{k_i}_{\rho}$, $i=1,2,\ldots,$ such that
$x_{i}\in
\partial\Omega_{\rho^*}$. Then $U(x_{i})=\rho^\ast$ for all $i$.
Since $\bar\Omega_{\rho^*}$
is compact, we may assume without loss of generality that
$\{x_{i}\}$ converges, say, to some $\bar x\in\bar\Omega_{\rho^*}$.
On one hand, we have $U(\bar x)=\rho^*$. But on the other hand,
since $\rho\geq U_{k_i} (x_{i})$ and $U_{k_i}\to U$
uniformly on $\bar\Omega_{\rho^*}$, taking limit $i\to\infty$
yields that
$\rho\ge U(\bar x)$. It follows that $\rho\ge\rho^*$, a
contradiction.
\end{pf}

%
%le3.2 #&#
\begin{Lemma} \label{Claim-2} $\Omega_\rho^k$ is a generalized Lipschitz
domain
for each $k\ge k(\rho)$.
\end{Lemma}

\begin{pf} We only consider the case $n>1$ because the case with $n=1$ is
trivial. Let $k\ge k(\rho)$ be fixed. We note by claim 1 that
$\partial\Omega^k_{\rho}$ is compact and contained in
$U_k^{-1}(\rho)$. Consider a point $x_0\in\partial\Omega^k_{\rho}$.
If $\nabla U_k(x_0)\ne0$, then the implicit function theorem
implies that, in a neighborhood
of $x_0$, $\partial\Omega_{\rho}^k$ is actually a $C^2$
hypersurface which coincides with
$U_k^{-1}(\rho)$. Let $\nabla U_k(x_0)= 0$. Then the Hessian
$D^2U_k(x_0)$ is nondegenerate because $U_k$ is a Morse function.
If $D^2U_k(x_0)$ is positive definite, then $x_0$ is a local minimal
point of $U_k$, and thus it cannot lie in $\bar\Omega_{\rho}^k$. If
$D^2U_k(x_0)$ is negative definite, then $x_0$ is a local maximal
point of $U_k$ and thus it must lie in the interior
$\Omega^k_\rho$.
%the interior of $\bar\Omega_{\rho}^k$ by definition of $\Omega^k_\rho$.
Hence,
$D^2U_k(x_0)$ must be a hyperbolic matrix.
Let $1\leq M<n$ be the number of positive eigenvalues of $D^2U_k(x_0)$.
Then by the Morse lemma \cite{GP}, there is a $C^2$ local change
of coordinates $v=(v_1,\ldots,v_n)=\Phi(x)$ in a neighborhood of
$x_0$ under which $\Phi(x_0)=0$ and
\[
U_k\bigl(\Phi^{-1}(v)\bigr)=\rho+ v_1^2+
\cdots+v_M^2-v_{M+1}^2-
\cdots-v_n^2.
\]
It follows that, near $x_0$, $\partial\Omega^k_\rho=
U_k^{-1}(\rho)$ is a union of two Lipschitz hypersurfaces
intersecting at $x_0$; each belongs to the boundary of a component
of $\Omega_{\rho}^k$. Since all nondegenerate critical points of
$U_k$ are isolated and $\partial\Omega_{\rho}^k$ is a compact set,
the number of critical points of $U_k$ on $\partial\Omega_{\rho}^k$
must be finite. Consequently, the number
of connected components of $\Omega_{\rho}^k$ which contain
nondegenerate critical points on their boundaries are finite. The
number of connected components of $\Omega_{\rho}^k$ which contain no
critical points on their boundaries is also finite, because each
such a component is separated from the rest of $\Omega_{\rho}^k$.
Thus, $\Omega^k_{\rho}$ is a generalized Lipschitz domain.
\end{pf}

\begin{pf*}{Proof of Theorem \ref{thmA}(a)} Let $\mu$ be a regular stationary
measure of \eqref{fp} with density $u\in W^{1,p}_{\mathrm{loc}}(\cal U)$.

For given $\rho_0\in(\rho_m,\rho_M)$, we consider a fixed
monotonically increasing function $\phi\in C^2(\R_+)$ satisfying
\[
\phi(t)=\cases{ 0, &\quad $\mbox{if $t\in[0,\rho_m]$;}$\vspace*{2pt}
\cr
t, &\quad $\mbox{if $t\in[\rho_0,+\infty)$.}$}
\]
We note that
$\phi^{\prime\prime}(t)=0$ for all $t\in[0,\rho_m]\cup
[\rho_0,+\infty)$.

Let $\rho\in(\rho_0,\rho_M)$ and $\rho^*\in(\rho,\rho_M)$. Since
$\phi\circ U_k\equiv\rho$ on $\partial\Omega_\rho^k$, using
Lemmas~\ref{Claim-1}, \ref{Claim-2}, we can apply Theorem \ref{l31}
with $F=\phi\circ U_k$, $\Omega=\Omega_{\rho^*}$, and
$\Omega'=\Omega_\rho^k$ for each $k\ge k(\rho)$ to obtain the
identity
\[
\di_{\Omega_\rho^k}\bigl(a^{ij} \,\partial^2_{ij}
\phi(U_k) + V^i \,\partial_i
\phi(U_k)\bigr)u \,\mathrm{ d}x=\di_{\partial\Omega_\rho^k}u a^{ij}
\,\partial_i \phi(U_k) \nu_j \,\mathrm{ d}s,
\]
that is,
%
%
%e3.1 #&#
\begin{eqnarray}\label{omega-1}
&&\di_{\Omega_\rho^k}\phi^\prime(U_k) ({\cal L}U_k)
u \,\mathrm{ d}x + \di_{\Omega_\rho^k} \phi^{\prime\prime}(U_k)
\bigl(a^{ij} \,\partial_i U_k\,
\partial_j U_k\bigr) u \,\mathrm{ d}x
\nonumber
\\[-8pt]
\\[-8pt]
\nonumber
&&\qquad=\di_{\partial\Omega_\rho^k}\phi^\prime(U_k)u a^{ij}\,
\partial_i U_k \nu_j \,\mathrm{ d}s=
\di_{\partial\Omega_\rho^k}u a^{ij}\,\partial_i U_k
\nu_j \,\mathrm{ d}s,
\end{eqnarray}
where $(\nu_j)$ denote the
unit outward normal vectors of $\partial\Omega^k_\rho$. For each
$k\ge k(\rho)$, if $\nabla U_k(x_0)\ne0$ at some $x_0\in
\partial\Omega_{\rho}^k$, then the implicit function theorem implies
that there is a neighborhood of $x_0$ on $\partial\Omega_{\rho}^k$
such that $\nu(x)=({\nabla U_k(x)})/({|\nabla U_k(x)|})$ within the
neighborhood. Thus,
\[
a^{ij}(x)\,\partial_i U_k (x)
\nu_j(x)\ge0,\qquad x\in\partial\Omega_{\rho}^k.
\]
It then follows from \eqref{omega-1} that
\[
\di_{\Omega_\rho^k}\phi^\prime(U_k) ({\cal L}U_k)
u \,\mathrm{ d}x+\di_{\Omega_\rho^k}\phi^{\prime\prime}(U_k)
\bigl(a^{ij} \,\partial_i U_k\,
\partial_j U_k\bigr) u \,\mathrm{ d}x \ge0,
\]
that is,
%
%
%e3.2 #&#
\begin{eqnarray}\label{omega-11}
&& \di_{\Omega_{\rho^*}\setminus
U^{-1}(\rho)}\chi_{\Omega_\rho^k}\phi^\prime(U_k) ({
\cal L}U_k)u \,\mathrm{ d}x +\di_{U^{-1}(\rho)\cap\Omega_\rho^k}
\phi^\prime(U_k)
({\cal L}U_k)u \,\mathrm{ d}x
\nonumber
\\
&&\qquad\ge- \di_{\Omega_{\rho^*}\setminus
U^{-1}(\rho)}\chi_{\Omega_\rho^k}\phi^{\prime\prime}(U_k)
\bigl(a^{ij} \,\partial_i U_k\,
\partial_j U_k\bigr) u \,\mathrm{ d}x\\
&&\qquad\quad{}-\di_{U^{-1}(\rho)\cap\Omega_\rho^k}\phi^{\prime
\prime}(U_k)
\bigl(a^{ij} \,\partial_i U_k\,
\partial_j U_k\bigr) u \,\mathrm{ d}x,\nonumber
\end{eqnarray}
where for any Borel set $E\subset\Omega_{\rho^*}$,
$\chi_{E}$ denotes the
indicator function of $E$
in
$\Omega_{\rho^*}$. Since $U$ is a Lyapunov function and $\rho\in
(\rho
_0,\rho_M)$,
we have
\begin{eqnarray*}
&&\di_{U^{-1}(\rho)\cap\Omega_\rho^k} \phi^\prime(U_k) ({\cal
L}U_k)u \,\mathrm{ d}x
\\
& &\qquad\le\di_{U^{-1}(\rho)\cap\Omega_\rho^k} \bigl|\phi^\prime(U_k){\cal
L}U_k-\phi^\prime(U){\cal L}U\bigr|u\,\mathrm{ d}x+\di_{U^{-1}(\rho)\cap
\Omega_\rho^k}({
\cal L}U)u\,\mathrm{ d}x
\\
&&\qquad \le\bigl(\bigl|\phi^\prime(U_k)-\phi^\prime(U)\bigr|_{C(\Omega_{\rho
^*})}|U|_{C^2(\Omega_{\rho^*})}+
\bigl|\phi^\prime(U_k)\bigr|_{C(\Omega_{\rho^*})}|U_k-U|_{C^2(\Omega_{\rho
^*})}
\bigr)
\\
&&\qquad\quad{} \times\di_{\Omega_{\rho^*}}\bigl(|A|+|V|\bigr)u \,\mathrm{ d}x -\gamma\mu
\bigl(U^{-1}(
\rho)\cap\Omega_\rho^k\bigr).
\end{eqnarray*}
It follows from the facts $u\in C(\bar\Omega_{\rho^*})$ and $U_k\to
U$ in $C^2(\bar\Omega_{\rho^*})$ that
%
%
%e3.3 #&#
\begin{equation}
\label{Lya-Uk} \limsup_{k\to\infty}\di_{U^{-1}(\rho)\cap\Omega
_\rho^k}
\phi^\prime(U_k) ({\cal L}U_k)u \,\mathrm{ d}x\leq0.
\end{equation}
Since $\phi^{\prime\prime}(\rho)=0$, we also have
%
%
%e3.4 #&#
\begin{eqnarray}\label{omega-111}
&&\lim_{k\to\infty}\di_{U^{-1}(\rho)\cap
\Omega_\rho^k}\bigl|\phi^{\prime\prime}(U_k)\bigr|
\bigl|a^{ij} \,\partial_i U_k \,\partial_j
U_k\bigr| u \,\mathrm{ d}x
\nonumber
\\
&&\qquad \le\lim_{k\to\infty}\di_{U^{-1}(\rho)}\bigl|\phi^{\prime\prime}(U_k)\bigr
|
\bigl|a^{ij} \,\partial_i U_k \,\partial_j
U_k\bigr| u \,\mathrm{ d}x
\\
& &\qquad\le\bigl|\phi^{\prime\prime}(\rho)\bigr||A|_{C(U^{-1}(\rho))}|\nabla
U|^2_{C(U^{-1}(\rho))}=0.\nonumber
\end{eqnarray}
Using the uniform convergence of $U_k\to U$
in $\Omega_{\rho^*}$,
it is easy to see that as
$ k\to\infty$,
%
%
%e3.5 #&#
\begin{equation}
\label{chi-k1} \chi_{\Omega_\rho^k}(x)\to\chi_{\Omega_\rho}(x),\qquad
x\in
\Omega_{\rho^*}\setminus U^{-1}(\rho).
\end{equation}
By taking limit $k\to\infty$ in \eqref{omega-11}
and using \eqref{Lya-Uk}--\eqref{chi-k1} and the dominated
convergence theorem, we now have
\begin{eqnarray*}
\di_{\Omega_{\rho}}\phi^\prime(U) ({\cal L}U)u \,\mathrm{ d}x&=&
\di_{\Omega_{\rho^*}\setminus
U^{-1}(\rho)}\chi_{\Omega_\rho}\phi^\prime(U) ({\cal L}U)u
\,\mathrm{ d}x
\\
&\ge&-\di_{\Omega_{\rho^*}\setminus U^{-1}(\rho)}\chi_{\Omega
_\rho} \phi^{\prime\prime}(U)
\bigl(a^{ij} \,\partial_i U \,\partial_j U\bigr) u
\,\mathrm{ d}x
\\
&=& -\di_{\Omega_\rho}\phi^{\prime\prime}(U) \bigl(a^{ij}
\,\partial_i U\, \partial_j U\bigr) u \,\mathrm{ d}x,
\end{eqnarray*}
which, by definition of $\phi$, is equivalent to
\[
\di_{\Omega_{\rho}\setminus\Omega_{\rho_m}}\phi^\prime(U)
({\cal L}U)u \,\mathrm{ d}x\ge-
\di_{\Omega_\rho\setminus
\Omega_{\rho_m}}\phi^{\prime\prime}(U) \bigl(a^{ij}\,
\partial_i U\, \partial_j U\bigr) u \,\mathrm{ d}x.
\]
Letting $\rho\to\rho_M$ in the above, we obtain
%
%
%e3.6 #&#
\begin{equation}
\label{28} \di_{{\cal U}\setminus\Omega_{\rho_m}}\phi^\prime(U)
({\cal L}U)u \,\mathrm{ d}x \ge-
\di_{{\cal U}\setminus\Omega_{\rho_m}} \phi^{\prime\prime}(U)
\bigl(a^{ij}\,
\partial_i U\, \partial_j U\bigr) u \,\mathrm{ d}x.
\end{equation}
We note that $\phi^\prime(t)\geq0$ and $\phi^\prime(t)=1$ as
$t\geq
\rho_0$. Using the fact that $U$ is a Lyapunov function, we clearly have
%
%
%e3.7 #&#
\begin{eqnarray}\label{est1}
\di_{{\cal U}\setminus\Omega_{\rho_m}}\phi^\prime(U) ({\cal L}U)u
\,\mathrm{ d}x &\leq&-\gamma
\di_{{\cal U}\setminus
\Omega_{\rho_m}}\phi^\prime(U)u \,\mathrm{ d}x
\nonumber
\\
&=&-\gamma\di_{{\cal U}\setminus\Omega_{\rho_0}} u \,\mathrm{ d}x-\gamma\di_{\Omega_{\rho_0} \setminus\Omega_{\rho_m}}
\phi^\prime(U)u \,\mathrm{ d}x
\nonumber
\\[-8pt]
\\[-8pt]
\nonumber
&\leq&-\gamma\di_{{\cal U}\setminus\Omega_{\rho_0}} u \,\mathrm{ d}x
\nonumber
\\
&=&-\gamma\mu({\cal U}\setminus\Omega_{\rho_0}).\nonumber
\end{eqnarray}
Denote $C_{\rho_m,\rho_0}=\max_{\rho_m\le\rho\le
\rho_0}|\phi^{\prime\prime}(\rho)|$. Then it is also clear that
%
%
%e3.8 #&#
\begin{eqnarray}\label{est2}
\di_{{\cal U}\setminus\Omega_{\rho_m}} \bigl|\phi^{\prime\prime}(U)\bigr|
\bigl(a^{ij}\,
\partial_i U\, \partial_j U\bigr) u \,\mathrm{ d}x&=&
\di_{\Omega_{\rho_0}\setminus\Omega_{\rho_m}} \bigl|\phi^{\prime
\prime}(U)\bigr| \bigl(a^{ij}\,
\partial_i U\, \partial_j U\bigr) u \,\mathrm{ d}x
\nonumber
\\
&\leq& C_{\rho_m, \rho_0} \Bigl(\sup_{\rho\in(\rho_m, \rho
_0)}H(\rho)\Bigr)\di
_{\Omega_{\rho_0}\setminus\Omega_{\rho_m}} u \,\mathrm{ d}x
\\
&=& C_{\rho_m,\rho_0} \Bigl(\sup_{\rho\in(\rho_m,
\rho_o)}H(\rho)\Bigr)\mu(
\Omega_{\rho_0}\setminus\Omega_{\rho_m}).\nonumber
\end{eqnarray}
The theorem now follows from (\ref{28})--\eqref{est2}.
\end{pf*}
%\qed

%s4 #&#
\section{Proof of Theorem \texorpdfstring{\protect\ref{thmA}(b)}{A(b)} and Theorem \texorpdfstring{\protect\ref{thmB}(a)}{B(a)}}\label{sec4}

Let $U$ be either a Lyapunov function or an anti-Lyapunov
function in $\cal U$ with respect to ${\cal L}$ with either
Lyapunov constant or anti-Lyapunov constant $\gamma$ and essential
lower, upper bound $\rho_m,\rho_M$, respectively, which satisfies
\eqref{regular-value-1}. Also let $H$ be as in \eqref{Hupdown} and
denote
$\Omega_\rho$ as the
$\rho$-sublevel set of $U$ for each $\rho\in[\rho_m,\rho_M)$. Let
$\mu$ be a regular stationary measure of \eqref{fp} with density
$u\in W^{1,p}_{\mathrm{loc}}(\cal U)$.

Consider the set $\cal I=\{\rho\in(\rho_m,\rho_M)\dvtx\nabla
U(x)\ne
0, x\in U^{-1}(\rho)\}$.
Then for each $\eta\in\cal I$,
$\Omega_\eta$ is a $C^2$ domain, whose boundary
$\partial\Omega_\eta$ coincides with $U^{-1}(\eta)$, and the
outward unit normal vector $\nu(x)$ of $\partial\Omega_\eta$ at each
$x$ is well defined and equals $({\nabla U(x)})/({|\nabla U(x)|})$.
Since $\cal I$ is open,
\[
\mathcal{I}=\bigcup_{1\le k<I}(a_k,b_k),
\]
where $I$ can be a positive integer or $+\infty$, and the intervals
$(a_k,b_k)$, $1\le k<I$, are pairwise disjoint.

\begin{pf*}{Proof of Theorem \ref{thmA}(b)} Let $\eta^*\in(\rho_m,\rho_M)\cap
\cal I$. For any $\eta\in(\rho_m, \eta^*)\cap\cal I$,
applications of Theorem \ref{l31} with $F=U$ on
$\Omega^{\prime}=\Omega_{\eta^*}$, $ \Omega_{\eta}$, respectively,
yield that
\begin{eqnarray*}
&&\di_{\partial\Omega_{\eta}}u a^{ij} \frac{\partial_i U\,
\partial_j U}{|\nabla U|} \,\mathrm{ d}s +
\di_{\Omega_{\eta^*}\setminus
\Omega_{\eta}}\bigl(a^{ij}\, \partial^2_{ij} U
+ V^i\, \partial_i U\bigr)u \,\mathrm{ d}x
\\
&&\qquad =\di_{\partial\Omega_{\eta^*}} ua^{ij}\frac{\partial_i U\,
\partial_j U}{|\nabla U|} \,\mathrm{ d}s.
\end{eqnarray*}
Since the right-hand side of the above is
nonnegative, applications of \eqref{Hupdown} to the first term of
the left-hand side of above and the definition of Lyapunov function
to the second term of the left-hand side of above yield that
%
%
%e4.1 #&#
\begin{equation}
\label{p120} \gamma\di_{\Omega_{\eta^*}\setminus\Omega_{\eta}}
u \,\mathrm{ d}x \leq H(\eta)\di_{\partial\Omega_\eta}
\frac{u}{|\nabla U|} \,\mathrm{ d}s,\qquad \eta\in\bigl[\rho_m,\eta^*\bigr)\cap
\cal I.
\end{equation}
Consider the function
\[
y(\eta)=\mu(\Omega_{\eta^*}\setminus\Omega_{\eta})=
\di_{\Omega_{\eta
^*}\setminus
\Omega_{\eta}} u \,\mathrm{ d}x,\qquad \eta\in\bigl(\rho_m,\eta^*\bigr
)\cap
\cal I. %
\]
By Theorem \ref{deri}, $y(\eta)$ is of the class $C^1$ on
$(\rho_m,\eta^*)\cap\cal I$ and
\[
y^{\prime}(\eta)= -\int_{\partial\Omega_\eta} \frac{u}{|\nabla
U|} \,\mathrm{ d}
s, \qquad\eta\in\bigl(\rho_m,\eta^*\bigr)\cap\cal I.
\]
Hence, by \eqref{p120},
%
%
%e4.2 #&#
\begin{equation}
\label{yeta} y^\prime(\eta)+\frac{\gamma} {H(\eta)} y(\eta)\le
0, \qquad\eta\in
\bigl(\rho_m,\eta^*\bigr)\cap\cal I.
\end{equation}
Let $1\leq k<I$ be
fixed. For any $ \widetilde{\eta}, \eta\in(a_k,b_k)$ with
$\widetilde{\eta}<\eta<\eta^*$, integrating \eqref{yeta} in the
interval $[\widetilde{\eta},\eta]$ yields that
\[
\mu(\Omega_{\eta^*}\setminus\Omega_{\eta})\le\mu(
\Omega_{\eta^*}\setminus\Omega_{\tilde\eta})\mathrm{ e}^{-\gamma
\int_{\tilde\eta}^\eta
{1}/{H(t)}\,\mathrm{ d}t}.
\]

In \eqref{yeta}, we have assumed without loss of generality
that $H$ is a positive function. If not, we can replace $H$ in \eqref
{p120} [hence in \eqref{yeta}] by
$H+\varepsilon$,
$0<\varepsilon\ll1$, so that the above estimate
holds with $H+\varepsilon$ in place of $H$. Since $y(\rho)$
is independent of $\varepsilon$, the estimate in fact holds for $H$
after taking $\varepsilon\to0$.

Since $\cal I$ is dense in $[\rho_m,\rho_M)$ by
\eqref{regular-value-1}, letting $\eta^*\to\rho_M$ in the above
yields that
%
%
%e4.3 #&#
\begin{equation}
\label{yeta1} \mu(\cal U\setminus\Omega_\eta)\le\mu(\cal
U\setminus
\Omega_{\tilde\eta})\mathrm{ e}^{-\gamma\int_{\tilde\eta}^\eta
{1}/{H(t)}\,\mathrm{ d}t}.
\end{equation}
By taking $\tilde\eta\to a_k$, or $\eta\to b_k$, and noting that
function $\mu(\cal U\setminus\Omega_t)$ is monotone in $t\in
[\rho_m,\rho_M]$, we see that \eqref{yeta1} in fact holds for all
$\tilde\eta, \eta\in[a_k,b_k]$ with
$\tilde\eta\le\eta$.

Next, let $\rho_*,\rho^*\in\cal I$ with $\rho_*<\rho^*$. We can
find $1\le\ell<I$ such that $\rho_*,\rho^*\in
\bigcup_{k=1}^\ell(a_k,b_k)$. Denote $I_\ell=\{ i\in
\{1,2,\ldots,\ell\}\dvtx(a_i,b_i)\cap[\rho_{*},\rho^*]\neq
\varnothing\}$
and $\tau=|I_\ell|$. Then $I_\ell=\{ i_1,i_2,\ldots,i_{\tau}\dvtx
b_{i_s}\le a_{i_{s+1}}, s=1,2,\ldots,\tau-1\}$, $\rho_*\in
(a_{i_1},b_{i_1})$, and $\rho^*\in(a_{i_{\tau}},b_{i_{\tau}})$. By
a recursive application of \eqref{yeta1} for
$k=i_1,i_2,\ldots,i_{\tau}$ respectively, we have
\begin{eqnarray*}
\mu(\cal U\setminus\Omega_{\rho^*})&\le&\mu(\cal U\setminus
\Omega_{a_{i_{\tau}}})\mathrm{ e}^{-\gamma\int_{a_{i_{\tau
}}}^{\rho^*}
{1}/{H(t)}\,\mathrm{ d}t}
\\
&\le&\mu(\cal U\setminus\Omega_{b_{i_{\tau-1}}})\mathrm{
e}^{-\gamma\int_{a_{i_{\tau}}}^{\rho^*}
{1}/{H(t)}\,\mathrm{ d}t}
\\
&\le&\mu(\cal U\setminus\Omega_{a_{i_{\tau-1}}})\mathrm{
e}^{-\gamma\int_{a_{i_{\tau-1}}}^{b_{i_{\tau-1}}}
{1}/{H(t)}\,\mathrm{ d}t}\cdot
\mathrm{ e}^{-\gamma\int_{a_{i_{\tau}}}^{\rho^*}
{1}/{H(t)}\,\mathrm{ d}t}
\\
&=&\mu(\cal U\setminus\Omega_{a_{i_{\tau-1}}})\mathrm{ e}^{-\gamma
\int_{[\rho_{*},\rho^*]\cap
\bigcup_{k=\tau-1}^{\tau} (a_{i_k},b_{i_k})}{1}/{H(t)}\,\mathrm{ d}t}
\\
&\le&\cdots\le\mu(\cal U\setminus\Omega_{a_{i_2}})\mathrm{
e}^{-\gamma\int_{[\rho_{*},\rho^*]\cap\bigcup
_{k=2}^{\tau} (a_{i_k},b_{i_k})} {1}/{H(t)}\,\mathrm{ d}t}
\\
& \le&\mu(\cal U\setminus\Omega_{b_{i_1}})\mathrm{ e}^{-\gamma\int
_{[\rho_{*},\rho^*]\cap\bigcup
_{k=2}^{\tau} (a_{i_k},b_{i_k})}
{1}/{H(t)}\,\mathrm{ d}t}
\\
& \le&\mu(\cal U\setminus\Omega_{\rho_*})\mathrm{ e}^{-\gamma\int
_{\rho_{*}}^{b_{i_1}} {1}/{H(t)}\,\mathrm{ d}t}\cdot
\mathrm{ e}^{-\gamma\int_{[\rho_{*},\rho^*]\cap
\bigcup_{k=2}^{\tau} (a_{i_k},b_{i_k})} {1}/{H(t)}\,\mathrm{ d}
t}
\\
&=&\mu(\cal U\setminus\Omega_{\rho_*})\mathrm{ e}^{-\gamma\int
_{[\rho_{*},\rho^*]\cap\bigcup
_{k=1}^{\tau} (a_{i_k},b_{i_k})} {1}/{H(t)}\,\mathrm{ d}t}
\\
&=&\mu(\cal U\setminus\Omega_{\rho_*})\mathrm{ e}^{-\gamma\int
_{[\rho_{*},\rho^*]\cap\bigcup
_{k=1}^{\ell} (a_{k},b_{k})} {1}/{H(t)}\,\mathrm{ d}t}
\\
&\le&\mathrm{ e}^{-\gamma\int_{[\rho_{*},\rho^*]\cap\bigcup
_{k=1}^{\ell} (a_{k},b_{k})} {1}/{H(t)}\,\mathrm{ d}t}.
\end{eqnarray*}
Let $\ell\to I$ in the above. Since $[\rho_{*},\rho^*]\cap\bigcup
_{k=1}^{\ell} (a_{k},b_{k})\to[\rho_{*},\rho^*]\cap\cal I$
which is a full Lebesgue measure subset of $[\rho_*,\rho^*]$, we
obtain
%
%
%e4.4 #&#
\begin{equation}
\label{yeta2} \mu(\cal U\setminus\Omega_{\rho^*})\le\mathrm{
e}^{-\gamma\int_{\rho_{*}}^{\rho^*} {1}/{H(t)}\,\mathrm{ d}t}.
\end{equation}

Now for any $\rho\in[\rho_m,\rho_M)$, we let $\rho_i^*,\rho^i_*$ be
sequences in $\cal I$ such that $\rho_i^*\nearrow\rho$ and $\rho^i_*
\searrow\rho_m$ as $i\to\infty$. Since \eqref{yeta2} holds with
$\rho^i_*,\rho_i^*$ in place of $\rho_*,\rho^*$ respectively for all
$i$, the proof is complete by taking $i\to\infty$.
\end{pf*}%\qed

\begin{pf*}{Proof of Theorem \ref{thmB}(a)}
% For any $\rho\in\cal I$, we recall
%that $\Omega_\rho$ is a $C^2$ domain and the outward unit normal
%vector $\nu(x)$ is well-defined and equals $\frac{\nabla
%U(x)}{|\nabla U(x)|}$ for each $x\in
%\partial\Omega_\rho$.
Let $\eta_*\in(\rho_m,\rho_M)\cap\cal
I$ and $\eta\in(\eta_*,\rho_M)\cap\cal I$ be arbitrarily chosen.
Applying Theorem \ref{l31} with $F=U$ on
$\Omega^{\prime}=\Omega_\eta$, $\Omega_{\eta_*}$, respectively, we
have
\begin{eqnarray*}
&&\di_{\partial\Omega_{\eta_*}}u a^{ij} \frac{\partial_i U\,
\partial_j U}{|\nabla U|} \,\mathrm{ d}s +
\di_{\Omega_\eta\setminus
\Omega_{\eta_*}}\bigl(a^{ij} \,\partial^2_{ij} U
+ V^i\, \partial_i U\bigr)u \,\mathrm{ d}x
\\
&&\qquad =\di_{\partial\Omega_\eta} ua^{ij}\frac{\partial_i U\,
\partial_j U}{|\nabla U|} \,\mathrm{ d}s.
\end{eqnarray*}
Since the first term in the left-hand side of above is
nonnegative, applications of the definition of anti-Lyapunov
function to the second term of the left-hand side of above and
\eqref{Hupdown} to the right-hand side of above yield that
%
%
%e4.5 #&#
\begin{equation}
\label{p12} \gamma\di_{\Omega_\eta\setminus\Omega_{\eta_*}} u
\,\mathrm{ d}x \leq H(\eta)\di_{\partial\Omega_\eta}
\frac{u}{|\nabla U|} \,\mathrm{ d}s.
\end{equation}
Consider the function
\[
y(\eta)=\mu(\Omega_{\eta}\setminus\Omega_{\eta_*})=
\di_{\Omega_\eta\setminus\Omega_{\eta_*}} u \,\mathrm{ d}x, \qquad\eta
\in(\eta_*,\rho_M). %
\]
Then by Theorem \ref{deri}, $y(\eta)$ is of class $C^1$ at each
$\eta\in\cal I\cap(\eta_*,\rho_M)$ with derivative
\[
y^{\prime}(\eta)= \int_{\partial\Omega_\eta} \frac{u}{|\nabla
U|} \,\mathrm{ d}s.
\]
Hence, \eqref{p12} yields that
%
%
%e4.6 #&#
\begin{equation}
\label{yeta-1-5.1} y^\prime(\eta)-\frac{\gamma}{H(\eta)} y(\eta
)\geq0,\qquad \eta\in(
\eta_*,\rho_M)\cap\cal I.
\end{equation}
Here, we have again assumed without loss of generality that $H$ is
a positive function, via the same reasoning as in the proof of
Theorem \ref{thmA}(b) above.

Fix $1\le k<I$. For any $\eta, \widetilde{\eta}\in(a_k,b_k)$ with
$\widetilde{\eta}<\eta$, we may assume that $\eta_*<
\widetilde{\eta}$. Integrating \eqref{yeta-1-5.1} in the interval
$[\widetilde{\eta},\eta]$ yields that
\[
\mu(\Omega_{\eta}\setminus\Omega_{\eta_*})\ge\mu(
\Omega_{\tilde\eta}\setminus\Omega_{\eta_*}) \mathrm{ e}^{\gamma
\int_{\tilde\eta}^\eta{1}/{H(t)}\,\mathrm{ d}t}.
\]
By
\eqref{regular-value-1}, $\cal I$ is dense in $[\rho_m,\rho_M)$.
Then by letting $\eta_*\searrow\rho_m$ in the above and noting that
$\lim_{\eta_*\searrow\rho_m }\mu(
\Omega_{\eta_*})=\mu(\Omega^*_{\rho_m})$, we have
%
%
%e4.7 #&#
\begin{equation}
\label{omegarhok-5.1}\mu\bigl(\Omega_{\eta}\setminus\Omega
^*_{\rho_m}
\bigr)\geq\mu\bigl(\Omega_{\tilde\eta}\setminus\Omega^*_{\rho
_m}\bigr)
\mathrm{ e}^{\gamma
\int_{\tilde\eta}^\eta{1}/{H(t)}\,\mathrm{ d}t},\qquad a_k< \tilde
\eta< \eta< b_k.
\end{equation}
We note that \eqref{omegarhok-5.1} also holds when $\tilde\eta=a_k$
or $ \eta=b_k$ by the monotonicity of the function $t\in
[\rho_m,\rho_M]\mapsto\mu(\Omega_t\setminus\Omega^*_{\rho_m})$.

Next, let $\rho_*,\rho^*\in\cal I$ with $\rho_*<\rho^*$. We fix
$1\le\ell<I$ such that $\rho_*,\rho^*\in\bigcup_{k=1}^\ell
(a_k,b_k)$. Denote $I_\ell=\{ i\in\{1,2,\ldots,\ell\}\dvtx
(a_i,b_i)\cap
[\rho_{*},\rho^*]\neq\varnothing\}$ and $\tau=|I_\ell|$. Then
$I_\ell=\{ i_1,i_2,\ldots,i_{\tau}\dvtx b_{i_s}\le a_{i_{s+1}},
s=1,2,\ldots,\tau-1\}$, $\rho_*\in(a_{i_1},b_{i_1})$, and
$\rho^*\in(a_{i_{\tau}},b_{i_{\tau}})$. In the case $\tau\geq2$,
by a recursive application of \eqref{omegarhok-5.1} for
$k=i_1,i_2,\ldots,i_{\tau}$ respectively, we have
\begin{eqnarray*}
\mu\bigl(\Omega_{\rho^*}\setminus\Omega^*_{\rho_m}\bigr)&\ge&
\mu
\bigl(\Omega_{a_{i_{\tau}}}\setminus\Omega^*_{\rho_m}\bigr)
\mathrm{ e}^{\gamma\int_{a_{i_{\tau}}}^{\rho^*}
{1}/{H(t)}\,\mathrm{ d}t}
\\
&\ge&\mu\bigl(\Omega_{b_{i_{\tau-1}}}\setminus\Omega^*_{\rho
_m}\bigr)
\mathrm{ e}^{\gamma\int_{a_{i_{\tau}}}^{\rho^*}
{1}/{H(t)}\,\mathrm{ d}t}
\\
&\ge&\mu\bigl(\Omega_{a_{i_{\tau-1}}}\setminus\Omega^*_{\rho
_m}\bigr)
\mathrm{ e}^{\gamma\int_{a_{i_{\tau-1}}}^{b_{i_{\tau-1}}}
{1}/{H(t)}\,\mathrm{ d}t}\cdot\mathrm{ e}^{\gamma\int
_{a_{i_{\tau}}}^{\rho^*}
{1}/{H(t)}\,\mathrm{ d}t}
\\
&=&\mu\bigl(\Omega_{a_{i_{\tau-1}}}\setminus\Omega^*_{\rho
_m}\bigr)
\mathrm{ e}^{\gamma\int_{[\rho_{*},\rho^*]\cap\bigcup
_{k=\tau-1}^{\tau} (a_{i_k},b_{i_k})}{1}/{H(t)}\,\mathrm{ d}t}
\\
&\ge&\cdots\ge\mu\bigl(\Omega_{a_{i_2}}\setminus\Omega^*_{\rho_m}
\bigr)\mathrm{ e}^{\gamma\int_{[\rho_{*},\rho^*]\cap\bigcup
_{k=2}^{\tau} (a_{i_k},b_{i_k})} {1}/{H(t)}\,\mathrm{ d}t}
\\
& \ge&\mu\bigl(\Omega_{b_{i_1}}\setminus\Omega^*_{\rho_m}\bigr)
\mathrm{ e}^{\gamma\int_{[\rho_{*},\rho^*]\cap\bigcup
_{k=2}^{\tau} (a_{i_k},b_{i_k})} {1}/{H(t)}\,\mathrm{ d}t}
\\
& \ge&\mu\bigl(\Omega_{\rho_*}\setminus\Omega^*_{\rho_m}\bigr)
\mathrm{ e}^{\gamma\int_{\rho_{*}}^{b_{i_1}}
{1}/{H(t)}\,\mathrm{ d}t}\cdot\mathrm{ e}^{\gamma\int_{[\rho
_{*},\rho^*]\cap
\bigcup_{k=2}^{\tau} (a_{i_k},b_{i_k})} {1}/{H(t)}\,\mathrm{ d}
t}
\\
&=&\mu\bigl(\Omega_{\rho_*}\setminus\Omega^*_{\rho_m}\bigr)
\mathrm{ e}^{\gamma\int_{[\rho_{*},\rho^*]\cap\bigcup
_{k=1}^{\tau} (a_{i_k},b_{i_k})} {1}/{H(t)}\,\mathrm{ d}t}
\\
&=&\mu\bigl(\Omega_{\rho_*}\setminus\Omega^*_{\rho_m}\bigr)
\mathrm{ e}^{\gamma\int_{[\rho_{*},\rho^*]\cap\bigcup
_{k=1}^{\ell} (a_{k},b_{k})} {1}/{H(t)}\,\mathrm{ d}t}.
\end{eqnarray*}
Let $\ell\to I$ in above. Since $[\rho_{*},\rho^*]\cap\cal I$ is of
full Lebesgue measure in $[\rho_{*},\rho^*]$, we have
%
%
%e4.8 #&#
\begin{equation}
\label{yeta2-5.1} \mu\bigl(\Omega_{\rho^*}\setminus\Omega
^*_{\rho_m}
\bigr)\ge\mu\bigl(\Omega_{\rho_*}\setminus\Omega^*_{\rho
_m}\bigr)
\mathrm{ e}^{\gamma\int_{\rho_{*}}^{\rho^*} {1}/{H(t)}\,\mathrm{ d}t}.
\end{equation}
In the case $\tau=1$, \eqref{yeta2-5.1} follows directly from
\eqref{omegarhok-5.1}.

Now for any $\rho_m<\rho_0<\rho<\rho_M$, we let $\rho^i_*,\rho_i^*$
be sequences in $\cal I$ such that $\rho_i^*\nearrow\rho$ and
$\rho^i_* \searrow\rho_0$ as $i\to\infty$. Since \eqref{yeta2-5.1}
holds with $\rho^i_*,\rho_i^*$ in place of $\rho_*,\rho^*$
respectively for all $i$, the proof is complete by taking
$i\to\infty$.
\end{pf*}

%s5 #&#
\section{Proof of Theorem \texorpdfstring{\protect\ref{thmA}(c)}{A(c)} and Theorem \texorpdfstring{\protect\ref{thmB}(b)}{B(b)}}\label{sec5}

Let $U$ be either a weak Lyapunov function or a weak anti-Lyapunov
function in $\cal U$ with respect to ${\cal L}$ with essential
lower, upper bound $\rho_m,\rho_M$, respectively. Also let $h,H$ be
as in \eqref{Hupdown} and denote
$\Omega_\rho$ as the
$\rho$-sublevel set of $U$ for each $\rho\in[\rho_m,\rho_M)$.

For each $\rho\in[\rho_m,\rho_M)$, since $h(\rho)>0$ in
\eqref{Hupdown},
$\nabla U(x)\ne0$ for all $x\in
U^{-1}(\rho)$ and $\Omega_\rho$ is a $C^2$ domain with
%
%
%e5.1 #&#
\begin{equation}
\label{Ome1}\partial\Omega_\rho=U^{-1}(\rho).
\end{equation}
Consider a regular stationary measure
$\mu$ of \eqref{fp} with density $u(x)\in W^{1,p}_{\mathrm{loc}}(\cal U)$.
Then by Theorem \ref{deri}, the function
\[
y(\rho)=\di_{\Omega_\rho\setminus
\Omega_{\rho_m}} u \,\mathrm{ d}x,\qquad \rho\in(\rho_m,
\rho_M)%
\]
is of the class $C^1$ and
\[
y^{\prime}(\rho)=\int_{\partial\Omega_{\rho}} \frac{u}{|\nabla
U|} \,\mathrm{ d}s,\qquad
\rho\in(\rho_m,\rho_M).
\]

For $t\in[\rho_m,\rho_M)$, consider
\[
H_*(t)=\cases{ %
\displaystyle\int_{\rho_m}^t
\int_{\rho_m}^\tau{1} {H(s)} \,\mathrm{ d}s\,\mathrm{ d}\tau, &\quad
$\mbox{in the case of Theorem \ref{thmA}(c)},$
\vspace*{2pt}\cr
\displaystyle\int_{\rho_m}^t\int_{\rho_m}^\tau
{1} {h(s)} \,\mathrm{ d}s\,\mathrm{ d}\tau, &\quad $\mbox{in the case of
Theorem \ref{thmB}(b)}.$}
\]
Since $H $ is positive and
continuous on $[\rho_m,\rho_M)$ in the case of Theorem \ref{thmA}(c) so is
$h$ in the case of
Theorem \ref{thmB}(b), $H_*$ is a $C^2$ function on $[\rho_m,\rho_M)$. We
extend $H_*$
to a $C^2$ function on $[0,\rho_M)$ and still denote it by $H_*$.

%
%le5.1 #&#
\begin{Lemma}\label{351} For each $\rho\in[\rho_m,\rho_M)$,
%
%
%e5.2 #&#
\begin{equation}
\label{eq351} \di_{\Omega_\rho\setminus
\Omega_{\rho_m}}\bigl(a^{ij} \,\partial^2_{ij}
F + V^i\,\partial_i F\bigr)u \,\mathrm{ d}x = H_*^\prime(
\rho)\di_{\partial\Omega_\rho}u a^{ij}\frac{\partial_i
U\,\partial_j U}{|\nabla U|} \,\mathrm{ d}s.
\end{equation}
\end{Lemma}

\begin{pf} Let $F=H_*\circ U$. For any $\rho\in
(\rho_m,\rho_M)$, we note that $F\in C^2(\bar\Omega_\rho)$ and
$F|_{\partial\Omega_\rho} \equiv H_*(\rho)$. We
apply Theorem \ref{l31} to $F$ with $\Omega^\prime$ being
$\Omega_\rho$, $\Omega_{\rho_m}$, respectively. By using the fact
that the unit outward normal vector ${\nu}(x)$ is well defined and
equals $\frac{\nabla U(x)}{|\nabla U(x)|}$ for any $ x\in
\partial\Omega_\rho\cup\partial\Omega_{\rho_m}$, we have
\begin{eqnarray*}
&&H_*^\prime(\rho_m)\di_{\partial\Omega_{\rho_m}}u a^{ij}
\frac{\partial_i U\,\partial_j U}{|\nabla U|} \,\mathrm{ d}s+\di
_{\Omega_\rho\setminus\Omega_{\rho_m}}\bigl(a^{ij}
\,\partial^2_{ij} F + V^i\,\partial_i F
\bigr)u \,\mathrm{ d}x
\\
&&\qquad =H_*^\prime(\rho)\di_{\partial\Omega_\rho}u a^{ij}
\frac{\partial_i
U\,\partial_j U}{|\nabla U|} \,\mathrm{ d}s.
\end{eqnarray*}
Since $H_*^\prime(\rho_m)=0$, the lemma holds.
\end{pf}

%Without loss of generality, we assume that $\rho_m\geq\rho_m^\prime$.
\begin{pf*}{Proof of Theorem \ref{thmA}(c)} To estimate the left-hand side of
(\ref{eq351}), we let $x\in\Omega_\rho\setminus\bar\Omega_{\rho_m}$
and denote $\rho^\prime=:U(x)$. Clearly, $\rho^\prime\in
(\rho_m,\rho)$. Then
$H_*^{\prime\prime}(\rho^\prime)=H^{-1}(\rho^\prime)$,
$H_*^{\prime}(\rho^\prime)\geq0$ and ${\cal L}U(x)=a^{ij}(x)
\,\partial^2_{ij} U(x) + V^i(x) \,\partial_i
U(x)\le0$. It then follows from~\eqref{Hupdown} that
\begin{eqnarray*}
a^{ij}(x) \,\partial^2_{ij} F(x) +
V^i(x)\,\partial_i F(x) &=&H_*^\prime\bigl(
\rho^\prime\bigr) \bigl(a^{ij}(x) \,\partial^2_{ij}
U(x) + V^i(x) \,\partial_i U(x)\bigr)
\\
&&{}+ H_*^{\prime\prime}\bigl(\rho^\prime\bigr)a^{ij}(x)
\partial_i U(x) \partial_j U(x) \\
&\leq&0 + \frac{1}{H(\rho^\prime)}
H\bigl(\rho^\prime\bigr)=1.
\end{eqnarray*}
Also by \eqref{Hupdown}, the right-hand side of (\ref{eq351}) simply
satisfies
\[
H_*^\prime(\rho)\di_{\partial\Omega_\rho} u a^{ij}
\frac{\partial_i
U\,
\partial_j U}{|\nabla U|} \,\mathrm{ d}s \geq\int_{\rho_m}^\rho
\frac{1}{H(s)}\,\mathrm{ d}s\, h(\rho) \int_{\partial\Omega_\rho
}\frac{
u}{|\nabla U|}
\,\mathrm{ d}s.
\]
Hence, by (\ref{eq351}),
\[
\di_{\Omega_\rho\setminus\Omega_{\rho_m}} u(x) \,\mathrm{ d}x\geq
\tilde{H}(\rho) \int_{\partial\Omega_\rho}
\frac{ u}{|\nabla U|} \,\mathrm{ d}s,\qquad \rho\in(\rho_m,\rho_M),
\]
that is,
\[
y^{\prime}(\rho)\le\frac{1}{\tilde
H(\rho)}y(\rho),\qquad \rho\in(\rho_m,
\rho_M).
\]

For any $\rho_0\in(\rho_m,\rho_M)$ and $\rho\in[\rho_0,\rho
_M)$, a
direct integration of the above inequality yields that
\[
y(\rho)\le y(\rho_0)\mathrm{ e}^{\int_{\rho_0}^{\rho}
{1}/{\tilde
H(r)}\,\mathrm{ d}r},\qquad \rho\in(
\rho_0,\rho_M).
\]
The proof is complete simply by taking limit $\rho\to\rho_M$ in
the above.
\end{pf*}

\begin{pf*}{Proof of Theorem \ref{thmB}(b)} To estimate the left-hand side of
(\ref{eq351}), we note that $H_*^{\prime}(t)\geq0$,
$H_*^{\prime\prime}(t)=h^{-1}(t)$ when $t>\rho_m$, and ${\cal
L}U=a^{ij}
\partial^2_{ij} U + V^i \partial_i
U\ge0$ in $\Omega_{\rho}\setminus\bar\Omega_{\rho_m}$. It then
follows from \eqref{Hupdown} that
\begin{eqnarray*}
a^{ij} \partial^2_{ij} F + V^i
\partial_i F &=&H_*^\prime(U) \bigl(a^{ij}
\partial^2_{ij} U + V^i \partial_i U
\bigr)+ H_*^{\prime\prime}(U)a^{ij}\partial_i U
\partial_j U
\\
&\geq& 0 + \frac{1}{h(U)} h(U)=1.
\end{eqnarray*}
Also by \eqref{Hupdown}, the right-hand side of (\ref{eq351}) simply
satisfies
\[
H_*^\prime(\rho)\di_{\partial\Omega_\rho} u a^{ij}
\frac{\partial_i
U\,
\partial_j U}{|\nabla U|} \,\mathrm{ d}s \leq\int_{\rho_m}^\rho
\frac{1}{h(s)}\,\mathrm{ d}s H(\rho) \int_{\partial\Omega_\rho
}\frac{
u}{|\nabla U|}
\,\mathrm{ d}s.
\]
Hence, (\ref{eq351}) becomes
\begin{eqnarray*}
\di_{\Omega_\rho\setminus\Omega_{\rho_m}}u(x) \,\mathrm{ d}x&\leq
& \int_{\rho_m}^\rho
\frac{1}{h(s)}\,\mathrm{ d}s H(\rho) \int_{\partial\Omega_\rho
}\frac{ u}{|\nabla U|}
\,\mathrm{ d}s
\\
&=& \tilde H(\rho) \int_{\partial\Omega_\rho}\frac{ u}{|\nabla
U|} \,\mathrm{ d}s,\qquad \rho
\in(\rho_m,\rho_M),
\end{eqnarray*}
that is,
\[
\frac{1}{\tilde H(\rho)} y(\rho)\le y^\prime(\rho),\qquad \rho\in(
\rho_m,\rho_M).
\]
For any $\rho_0\in(\rho_m,\rho_M)$, let $\rho\in[\rho_0,\rho_M)$
be fixed. The proof is complete by
a direct integration of the above in the interval $[\rho_0,\rho]$.
\end{pf*}
%\qed

\section*{Acknowledgement} We would like to thank the referees
for valuable comments and also for letting us know the references
\cite{ABG,BKS1,BRS2,BRS3}.

% imsref loaded by akundreckaite, 2014-03-06 10:48:41
%

% zodis "Acknowledgments" paliekamas pagal autoriu
%\section*{Acknowledgments}

%\begin{supplement}[id=suppA]
%\sname{Supplement A}
%\stitle{}
%\slink[doi]{10.1214/00-AOPXXXXSUPP} %[doi,text={...}] - jei reikia
%suskaldyti doi
%\sdatatype{.pdf}
%\sfilename{aopXXXX\_supp.pdf}
%\sdescription{}
%\end{supplement}

%\begin{thebibliography}{99}
%\bibitem[\protect\citeauthoryear{}{}]{r1}
%\bibitem{r1}
%\end{thebibliography}

\printaddresses


\begin{thebibliography}{25}
% pybtex-1.00. Style name=ims, version=2.7, label_style=nolabel,
%sorting_style=complex, cfg=None, language=None.

%b1 ###
%b1 #&#
\bibitem{ABR}
%
\begin{barticle}[mr]
\bauthor{\bsnm{Albeverio},~\bfnm{S.}\binits{S.}},
\bauthor{\bsnm{Bogachev},~\bfnm{V.}\binits{V.}} \AND
\bauthor{\bsnm{R{\"o}ckner},~\bfnm{M.}\binits{M.}}
(\byear{1999}).
\btitle{On uniqueness of invariant measures for finite- and
infinite-dimensional diffusions}.
\bjournal{Comm. Pure Appl. Math.}
\bvolume{52}
\bpages{325--362}.
\bid
{doi={10.1002/(SICI)1097-0312(199903)52:3<325::AID-CPA2>3.0.CO;2-V},
issn={0010-3640}, mr={1656067}}
\end{barticle}
%
\bptok{imsref}%
% NOT OUTPUTED:
% issn = 0010-3640
% url =
%%http://dx.doi.org/10.1002/(SICI)1097-0312(199903)52:3<325::AID-CPA2>3.0.CO;2-V
% number = 3
% coden = CPAMA
% fjournal = Communications on Pure and Applied Mathematics
\endbibitem

%b2 ###
%b2 #&#
\bibitem{ABG}
%
\begin{bbook}[mr]
\bauthor{\bsnm{Arapostathis},~\bfnm{Ari}\binits{A.}},
\bauthor{\bsnm{Borkar},~\bfnm{Vivek~S.}\binits{V.~S.}} \AND
\bauthor{\bsnm{Ghosh},~\bfnm{Mrinal~K.}\binits{M.~K.}}
(\byear{2012}).
\btitle{Ergodic Control of Diffusion Processes}.
\bseries{Encyclopedia of Mathematics and Its Applications}
\bvolume{143}.
\bpublisher{Cambridge Univ. Press},
\blocation{Cambridge}.
\bid{mr={2884272}}
\end{bbook}
%
\bptok{imsref}%
% NOT OUTPUTED:
% isbn = 978-0-521-76840-5
% fpage = xvi+323
\endbibitem

%b3 ###
%b3 #&#
\bibitem{Ben}
%
\begin{bbook}[mr]
\bauthor{\bsnm{Bensoussan},~\bfnm{Alain}\binits{A.}}
(\byear{1988}).
\btitle{Perturbation Methods in Optimal Control}.
%\bseries{Wiley/Gauthier-Villars Series in Modern Applied Mathematics}.
\bpublisher{Wiley},
\blocation{Chichester}.
%\bnote{Translated from the French by C. Tomson}.
\bid{mr={0949208}}
\end{bbook}
%
\bptok{imsref}%
% NOT OUTPUTED:
% isbn = 0-471-91994-2
% fpage = xiv+573
\endbibitem

%b4 ###
%b4 #&#
\bibitem{NEW}
%
\begin{barticle}[mr]
\bauthor{\bsnm{Bhattacharya},~\bfnm{R.~N.}\binits{R.~N.}}
(\byear{1978}).
\btitle{Criteria for recurrence and existence of invariant measures for
multidimensional diffusions}.
\bjournal{Ann. Probab.}
\bvolume{6}
\bpages{541--553}.
\bid{issn={0091-1798}, mr={0494525}}
\end{barticle}
%
\bptok{imsref}%
% NOT OUTPUTED:
% issn = 0091-1798
% number = 4
% fjournal = The Annals of Probability
\endbibitem

%b5 ###
%b5 #&#
\bibitem{Boga}
%
\begin{bbook}[mr]
\bauthor{\bsnm{Bogachev},~\bfnm{V.~I.}\binits{V.~I.}}
(\byear{2007}).
\btitle{Measure Theory. {V}ol. {I}, {II}}.
\bpublisher{Springer},
\blocation{Berlin}.
\bid{doi={10.1007/978-3-540-34514-5}, mr={2267655}}
\end{bbook}
%
\bptok{imsref}%
% NOT OUTPUTED:
% isbn = 978-3-540-34513-8; 3-540-34513-2
% url = http://dx.doi.org/10.1007/978-3-540-34514-5
% fpage = Vol. I: xviii+500 pp., Vol. II: xiv+575
\endbibitem

%b6 ###
%b6 #&#
\bibitem{BDR}
%
\begin{barticle}[mr]
\bauthor{\bsnm{Bogachev},~\bfnm{V.~I.}\binits{V.~I.}},
\bauthor{\bsnm{Da Prato},~\bfnm{G.}\binits{G.}} \AND
\bauthor{\bsnm{R{\"o}ckner},~\bfnm{M.}\binits{M.}}
(\byear{2008}).
\btitle{On parabolic equations for measures}.
\bjournal{Comm. Partial Differential Equations}
\bvolume{33}
\bpages{397--418}.
\bid{doi={10.1080/03605300701382415}, issn={0360-5302}, mr={2398235}}
\end{barticle}
%
\bptok{imsref}%
% NOT OUTPUTED:
% issn = 0360-5302
% url = http://dx.doi.org/10.1080/03605300701382415
% number = 1-3
% coden = CPDIDZ
% fjournal = Communications in Partial Differential Equations
\endbibitem

%b7 ###
%b7 #&#
\bibitem{BKS1}
%
\begin{barticle}[mr]
\bauthor{\bsnm{Bogachev},~\bfnm{V.~I.}\binits{V.~I.}},
\bauthor{\bsnm{Kirillov},~\bfnm{A.~I.}\binits{A.~I.}} \AND
\bauthor{\bsnm{Shaposhnikov},~\bfnm{S.~V.}\binits{S.~V.}}
(\byear{2012}).
\btitle{Integrable solutions of the stationary {K}olmogorov equation}.
\bjournal{Dokl. Math.}
\bvolume{85}
\bpages{309--314}.
%\bid{doi={10.1134/S1064562412030027}, issn={0869-5652}, mr={2985904}}
\end{barticle}
%
\bptok{imsref}%
% NOT OUTPUTED:
% issn = 0869-5652
% url = http://dx.doi.org/10.1134/S1064562412030027
% number = 1
% fjournal = Rossi\u\i skaya Akademiya Nauk. Doklady Akademii Nauk
\endbibitem

%b8 ###
%b8 #&#
\bibitem{BKR09}
%
\begin{barticle}[mr]
\bauthor{\bsnm{Bogachev},~\bfnm{V.~I.}\binits{V.~I.}},
\bauthor{\bsnm{Krylov},~\bfnm{N.~V.}\binits{N.~V.}} \AND
\bauthor{\bsnm{R{\"o}kner},~\bfnm{M.}\binits{M.}}
(\byear{2009}).
\btitle{Elliptic and parabolic equations for measures}.
\bjournal{Russian Math. Surveys}
\bvolume{64}
\bpages{973--1078}.
%\bid{doi={10.1070/RM2009v064n06ABEH004652}, issn={0042-1316},
%mr={2640966}}
\end{barticle}
%
\bptok{imsref}%
% NOT OUTPUTED:
% issn = 0042-1316
% url = http://dx.doi.org/10.1070/RM2009v064n06ABEH004652
% number = 6(390)
% fjournal = Rossi\u\i skaya Akademiya Nauk. Moskovskoe Matematicheskoe
%Obshchestvo. Uspekhi Matematicheskikh Nauk
\endbibitem

%b9 ###
%b9 #&#
\bibitem{BKR}
%
\begin{barticle}[mr]
\bauthor{\bsnm{Bogachev},~\bfnm{V.~I.}\binits{V.~I.}},
\bauthor{\bsnm{Krylov},~\bfnm{N.~V.}\binits{N.~V.}} \AND
\bauthor{\bsnm{R{\"o}ckner},~\bfnm{M.}\binits{M.}}
(\byear{2001}).
\btitle{On regularity of transition probabilities and invariant
measures of singular diffusions under minimal conditions}.
\bjournal{Comm. Partial Differential Equations}
\bvolume{26}
\bpages{2037--2080}.
\bid{doi={10.1081/PDE-100107815}, issn={0360-5302}, mr={1876411}}
\end{barticle}
%
\bptok{imsref}%
% NOT OUTPUTED:
% issn = 0360-5302
% url = http://dx.doi.org/10.1081/PDE-100107815
% number = 11-12
% coden = CPDIDZ
% fjournal = Communications in Partial Differential Equations
\endbibitem

%b10 ###
%b10 #&#
\bibitem{BR01}
%
\begin{barticle}[mr]
\bauthor{\bsnm{Bogachev},~\bfnm{V.~I.}\binits{V.~I.}} \AND
\bauthor{\bsnm{R{\"o}kner},~\bfnm{M.}\binits{M.}}
(\byear{2001}).
\btitle{A generalization of Has'minski{\u\i}'s theorem on the existence
of invariant measures for locally integrable drifts}.
\bjournal{Theory Probab. Appl.}
\bvolume{45}
\bpages{363--378}.
%\bid{doi={10.1137/S0040585X97978348}, issn={0040-361X}, mr={1967783}}
\bptnote{check year}%
\end{barticle}
%
\bptok{imsref}%
% NOT OUTPUTED:
% issn = 0040-361X
% url = http://dx.doi.org/10.1137/S0040585X97978348
% number = 3
% fjournal = Rossi\u\i skaya Akademiya Nauk. Teoriya Veroyatnoste\u\i\
%i ee Primeneniya
\endbibitem

%b11 ###
%b11 #&#
\bibitem{BRS2}
%
\begin{barticle}[mr]
\bauthor{\bsnm{Bogachev},~\bfnm{V.~I.}\binits{V.~I.}},
\bauthor{\bsnm{R{\"o}kner},~\bfnm{M.}\binits{M.}} \AND
\bauthor{\bsnm{Shaposhnikov},~\bfnm{S.~V.}\binits{S.~V.}}
(\byear{2012}).
\btitle{On positive and probability solutions of the stationary
{F}okker--{P}lanck--{K}olmogorov equation}.
\bjournal{Dokl. Math.}
\bvolume{85}
\bpages{350--354}.
%\bid{doi={10.1134/S1064562412030143}, issn={0869-5652}, mr={2986416}}
\end{barticle}
%
\bptok{imsref}%
% NOT OUTPUTED:
% issn = 0869-5652
% url = http://dx.doi.org/10.1134/S1064562412030143
% number = 3
% fjournal = Rossi\u\i skaya Akademiya Nauk. Doklady Akademii Nauk
\endbibitem

%b12 ###
%b12 #&#
\bibitem{BRS}
%
\begin{barticle}[mr]
\bauthor{\bsnm{Bogachev},~\bfnm{V.~I.}\binits{V.~I.}},
\bauthor{\bsnm{R{\"o}kner},~\bfnm{M.}\binits{M.}} \AND
\bauthor{\bsnm{Shtannat},~\bfnm{V.}\binits{V.}}
(\byear{2002}).
\btitle{Uniqueness of solutions of elliptic equations and uniqueness of
invariant measures of diffusions}.
\bjournal{Mat. Sb.}
\bvolume{193}
\bpages{3--36}.
\bid{doi={10.1070/SM2002v193n07ABEH000665}, issn={0368-8666}, mr={1936848}}
\end{barticle}
%
\bptok{imsref}%
% NOT OUTPUTED:
% issn = 0368-8666
% url = http://dx.doi.org/10.1070/SM2002v193n07ABEH000665
% number = 7
% fjournal = Rossi\u\i skaya Akademiya Nauk. Matematicheski\u\i\ Sbornik
\endbibitem

%b13 ###
%b13 #&#
\bibitem{BR02}
%
\begin{bincollection}[mr]
\bauthor{\bsnm{Bogachev},~\bfnm{Vladimir~I.}\binits{V.~I.}} \AND
\bauthor{\bsnm{R{\"o}ckner},~\bfnm{Michael}\binits{M.}}
(\byear{2002}).
\btitle{Invariant measures of diffusion processes: Regularity,
existence, and uniqueness problems}.
In \bbooktitle{Stochastic Partial Differential Equations and
Applications ({T}rento, 2002)}.
\bseries{Lecture Notes in Pure and Applied Mathematics}
\bvolume{227}
\bpages{69--87}.
\bpublisher{Dekker},
\blocation{New York}.
\bid{mr={1919503}}
\end{bincollection}
%
\bptok{imsref}%
\endbibitem

%b14 ###
%b14 #&#
\bibitem{BRS3}
%
\begin{barticle}[mr]
\bauthor{\bsnm{Bogachev},~\bfnm{V.~I.}\binits{V.~I.}},
\bauthor{\bsnm{R{\"o}ckner},~\bfnm{M.}\binits{M.}} \AND
\bauthor{\bsnm{Shaposhnikov},~\bfnm{S.~V.}\binits{S.~V.}}
(\byear{2011}).
\btitle{On uniqueness problems related to elliptic equations for measures}.
\bjournal{J. Math. Sci. (N. Y.)}
\bvolume{176}
\bpages{759--773}.
%\bnote{Problems in mathematical analysis. No. 58}.
\bid{doi={10.1007/s10958-011-0434-3}, issn={1072-3374}, mr={2838973}}
\end{barticle}
%
\bptok{imsref}%
% NOT OUTPUTED:
% issn = 1072-3374
% url = http://dx.doi.org/10.1007/s10958-011-0434-3
% number = 6
% coden = JMTSEW
% fjournal = Journal of Mathematical Sciences (New York)
\endbibitem

%b15 ###
%b15 #&#
\bibitem{BRS1}
%
\begin{bincollection}[mr]
\bauthor{\bsnm{Bogachev},~\bfnm{Vladimir~I.}\binits{V.~I.}},
\bauthor{\bsnm{R{\"o}ckner},~\bfnm{Michael}\binits{M.}} \AND
\bauthor{\bsnm{Stannat},~\bfnm{Wilhelm}\binits{W.}}
(\byear{2000}).
\btitle{Uniqueness of invariant measures and essential
{$m$}-dissipativity of diffusion operators on {$L^1$}}.
In \bbooktitle{Infinite Dimensional Stochastic Analysis ({A}msterdam, 1999)}.
\bseries{Verh. Afd. Natuurkd. 1. Reeks. K. Ned. Akad. Wet.}
\bvolume{52}
\bpages{39--54}.
\bpublisher{R. Neth. Acad. Arts Sci.},
\blocation{Amsterdam}.
\bid{mr={1831410}}
\end{bincollection}
%
\bptok{imsref}%
\endbibitem

%b16 ###
%b16 #&#
\bibitem{GP}
%
\begin{bbook}[mr]
\bauthor{\bsnm{Guillemin},~\bfnm{Victor}\binits{V.}} \AND
\bauthor{\bsnm{Pollack},~\bfnm{Alan}\binits{A.}}
(\byear{1974}).
\btitle{Differential Topology}.
\bpublisher{Prentice-Hall},
\blocation{Englewood Cliffs, NJ.}
\bid{mr={0348781}}
\end{bbook}
%
\bptok{imsref}%
% NOT OUTPUTED:
% fpage = xvi+222
\endbibitem

%b17 ###
%b17 #&#
\bibitem{Ha2}
%
\begin{barticle}[mr]
\bauthor{\bsnm{Has'minski{\u\i}},~\bfnm{R.~Z.}\binits{R.~Z.}}
(\byear{1960}).
\btitle{Ergodic properties of recurrent diffusion processes and
stabilization of the solution of the {C}auchy problem for parabolic equations}.
\bjournal{Theory Probab. Appl.}
\bvolume{5}
\bpages{179--196}.
%\bid{issn={0040-361X}, mr={0133871}}
\end{barticle}
%
\bptok{imsref}%
% NOT OUTPUTED:
% issn = 0040-361x
% fjournal = Akademija Nauk SSSR. Teorija Verojatnoste\u\i\ i ee
%Primenenija
\endbibitem

%b18 ###
%b18 #&#
\bibitem{Ha}
%
\begin{bbook}[mr]
\bauthor{\bsnm{Has'minski{\u\i}},~\bfnm{R.~Z.}\binits{R.~Z.}}
(\byear{1980}).
\btitle{Stochastic Stability of Differential Equations}.
\bseries{Monographs and Textbooks on Mechanics of Solids and Fluids:
Mechanics and Analysis}
\bvolume{7}.
\bpublisher{Sijthoff \& Noordhoff},
\blocation{Alphen aan den Rijn}.
%\bnote{Translated from the Russian by D. Louvish}.
\bid{mr={0600653}}
\end{bbook}
%
\bptok{imsref}%
% NOT OUTPUTED:
% isbn = 90-286-0100-7
% fpage = xvi+344
\endbibitem

%b19 ###
%b19 #&#
\bibitem{JM1}
%
\begin{bmisc}[auto:STB|2014/02/12|14:17:21]
\bauthor{\bsnm{Huang},~\bfnm{W.}\binits{W.}},
\bauthor{\bsnm{Ji},~\bfnm{M.}\binits{M.}},
\bauthor{\bsnm{Liu},~\bfnm{Z.}\binits{Z.}} \AND
\bauthor{\bsnm{Yi},~\bfnm{Y.}\binits{Y.}}
(\byear{2013}).
\bhowpublished{Steady states of Fokker--Planck equations, Parts I--III. Submitted}.
\end{bmisc}
%
\bptok{imsref}%
% NOT OUTPUTED:
% sortkey = Huang
\endbibitem

%b20 ###
%b20 #&#
\bibitem{JM3}
%
\begin{bmisc}[auto:STB|2014/02/12|14:17:21]
\bauthor{\bsnm{Huang},~\bfnm{W.}\binits{W.}},
\bauthor{\bsnm{Ji},~\bfnm{M.}\binits{M.}},
\bauthor{\bsnm{Liu},~\bfnm{Z.}\binits{Z.}} \AND
\bauthor{\bsnm{Yi},~\bfnm{Y.}\binits{Y.}}
(\byear{2013}).
\bhowpublished{Concentration and limit behaviors of stationary
measures. Preprint}.
\end{bmisc}
%
\bptok{imsref}%
% NOT OUTPUTED:
% sortkey = Huang
\endbibitem

%b21 ###
%b21 #&#
\bibitem{JM4}
%
\begin{bmisc}[auto:STB|2014/02/12|14:17:21]
\bauthor{\bsnm{Huang},~\bfnm{W.}\binits{W.}},
\bauthor{\bsnm{Ji},~\bfnm{M.}\binits{M.}},
\bauthor{\bsnm{Liu},~\bfnm{Z.}\binits{Z.}} \AND
\bauthor{\bsnm{Yi},~\bfnm{Y.}\binits{Y.}}
(\byear{2014}).
\bhowpublished{Convergence of Gibbs measures. Submitted}.
\end{bmisc}
%
\bptok{imsref}%
% NOT OUTPUTED:
% sortkey = Huang
\endbibitem

%b22 ###
%b22 #&#
\bibitem{Sko}
%
\begin{bbook}[auto:STB|2014/02/12|14:17:21]
\bauthor{\bsnm{Skorohod},~\bfnm{A.~V.}\binits{A.~V.}}
(\byear{1989}).
\btitle{Asymptotic Methods in the Theory of Stochastic Differential Equations}.
\bpublisher{Amer. Math. Soc.},
\blocation{Providence, RI}.
\bid{mr={1020057}}
\end{bbook}
%
\bptok{imsref}%
\endbibitem

%b23 ###
%b23 #&#
\bibitem{Ver87}
%
\begin{barticle}[auto]
\bauthor{\bsnm{Veretennikov},~\bfnm{A.~Y.}\binits{A.~Y.}}
(\byear{1987}).
\btitle{Bounds for the mixing rate in the theory of
stochastic equations}.
\bjournal{Theory Probab. Appl.}
\bvolume{32}
\bpages{273--281}.
%\bid{issn={0040-361X}, mr={0902757}}
\end{barticle}
%
\bptok{imsref}%
% NOT OUTPUTED:
% issn = 0040-361X
% number = 2
% fjournal = Akademiya Nauk SSSR. Teoriya Veroyatnoste\u\i\ i ee
%Primeneniya
\endbibitem

%b24 ###
%b24 #&#
\bibitem{Ver97}
%
\begin{barticle}[mr]
\bauthor{\bsnm{Veretennikov},~\bfnm{A.~Y.}\binits{A.~Y.}}
(\byear{1997}).
\btitle{On polynomial mixing bounds for stochastic differential equations}.
\bjournal{Stochastic Process. Appl.}
\bvolume{70}
\bpages{115--127}.
\bid{doi={10.1016/S0304-4149(97)00056-2}, issn={0304-4149}, mr={1472961}}
\end{barticle}
%
\bptok{imsref}%
% NOT OUTPUTED:
% issn = 0304-4149
% url = http://dx.doi.org/10.1016/S0304-4149(97)00056-2
% number = 1
% coden = STOPB7
% fjournal = Stochastic Processes and their Applications
\endbibitem

%b25 ###
%b25 #&#
\bibitem{Ver99}
%
\begin{barticle}[mr]
\bauthor{\bsnm{Veretennikov},~\bfnm{A.~Y.}\binits{A.~Y.}}
(\byear{1999}).
\btitle{On polynomial mixing and the rate of convergence for stochastic
differential and difference equations}.
\bjournal{Theory Probab. Appl.}
\bvolume{44}
\bpages{361--374}.
%\bid{doi={10.1137/S0040585X97977550}, issn={0040-361X}, mr={1751475}}
\end{barticle}
%
\bptok{imsref}%
% NOT OUTPUTED:
% issn = 0040-361X
% url = http://dx.doi.org/10.1137/S0040585X97977550
% number = 2
% fjournal = Rossi\u\i skaya Akademiya Nauk. Teoriya Veroyatnoste\u\i\
%i ee Primeneniya
\endbibitem

\end{thebibliography}
\end{document}